\documentclass[preprint,10pt,1p]{elsarticle}

\usepackage{amssymb}
\usepackage[T1]{fontenc} 
\usepackage{fourier} 
\usepackage[english]{babel} 
\usepackage{amsmath,amsfonts,amsthm} 
\usepackage{appendix}

\usepackage{amssymb}
\usepackage[all]{xy}
\usepackage{enumerate}
\usepackage{tikz}
\usepackage{mathrsfs}
\usepackage[english]{babel}
\usepackage{multirow}
\usepackage{float}
\usepackage{enumitem}
\usepackage{graphicx} 

\newtheorem{innercustomthm}{Theorem}
\newenvironment{customthm}[1]
  {\renewcommand\theinnercustomthm{#1}\innercustomthm}
  {\endinnercustomthm}

\newtheorem{innercustomlemma}{Lemma}
\newenvironment{customlemma}[1]
  {\renewcommand\theinnercustomlemma{#1}\innercustomlemma}
  {\endinnercustomlemma}

\newtheorem{innercustomcorollary}{Corollary}
\newenvironment{customcorollary}[1]
  {\renewcommand\theinnercustomcorollary{#1}\innercustomcorollary}
  {\endinnercustomcorollary}

\newtheorem{innercustomprop}{Proposition}
\newenvironment{customprop}[1]
  {\renewcommand\theinnercustomprop{#1}\innercustomprop}
  {\endinnercustomprop}

\theoremstyle{definition}
\newtheorem{innercustomexample}{Example}
\newenvironment{customexample}[1]
  {\renewcommand\theinnercustomexample{#1}\innercustomexample}
  {\endinnercustomexample}

\newtheorem{innercustomremark}{Remark}
\newenvironment{customremark}[1]
  {\renewcommand\theinnercustomremark{#1}\innercustomremark}
  {\endinnercustomremark}

\newtheorem{innercustomconjecture}{Conjecture}
\newenvironment{customconjecture}[1]
{\renewcommand\theinnercustomconjecture{#1}\innercustomconjecture}
{\endinnercustomconjecture}

\newtheorem{innercustomquestion}{Question}
\newenvironment{customquestion}[1]
{\renewcommand\theinnercustomquestion{#1}\innercustomquestion}
{\endinnercustomquestion}

  \makeatletter
    \def\ps@pprintTitle{%
      \let\@oddhead\@empty
      \let\@evenhead\@empty
      \let\@oddfoot\@empty
      \let\@evenfoot\@oddfoot
    }
    \makeatother

\newcommand{\indicator}{{\mathbb{1}}}

\newcommand{\floor}[1]{\left\lfloor#1\right\rfloor}

\def\nbe{\mathtt{ne}}
\def\entropy{\mathtt{h}}
\def\lentropy{\underline{\mathtt{h}}}
\def\br{\mathtt{br}}
\def\rentropy{\overline{\mathtt{h}}}

\begin{document}

\begin{frontmatter}



\title{Limiting behavior of 3-color excitable media on arbitrary graphs}


\author{Janko Gravner\footnote{Department of Mathematics, University of California, 
Davis, CA 95616, USA, gravner@math.ucdavis.edu}}

\author{Hanbaek Lyu\footnote{Department of Mathematics, The Ohio State University, Columbus, OH 43210, yu.1242@osu.edu}}

\author{and David Sivakoff\footnote{Departments of Statistics and Mathematics, The Ohio State University, Columbus, OH 43210, \\ \text{\qquad dsivakoff@stat.osu.edu}}}

\newtheorem{problem}{Problem}
\newtheorem{proposition}{Proposition}[section]
\newtheorem{lemma}{Lemma}[section]
\newtheorem{corollary}{Corollary}[section]
\newtheorem{definition}{Definition}[section]
\newtheorem{theorem}{Theorem}[section]
\newtheorem{remark}{Remark}[section]
\newtheorem{ex}{Example}
\newtheorem{exercise}{Exercise}
\newtheorem{note}{Note}

\def\bibsection{\section*{References}}

\begin{abstract}
	Fix a simple graph $G=(V,E)$ and choose a random initial 3-coloring of vertices drawn from a uniform product measure. The 3-color cycle cellular automaton is a process in which at each discrete time step in parallel, every vertex with color $i$ advances to the successor color $(i+1)$ mod 3 if in contact with a neighbor with the successor color, and otherwise retains the same color. In the Greenberg-Hastings Model, the same update rule applies only to color 0, while other two colors automatically advance. The limiting behavior of these processes has been studied mainly on the integer lattices. In this paper, we introduce a monotone comparison process defined on the universal covering space of the underlying graph, and characterize the limiting behavior of these processes on arbitrary connected graphs. In particular, we establish a phase transition on the Erd\"os-R\'enyi random graph. On infinite trees, we connect the rate of color change to the cloud speed of an associated tree-indexed walk. We give estimates of the cloud speed by generalizing known results to trees with leaves. 
\end{abstract}


\end{frontmatter}


\section{Introduction}

An excitable medium is a network of coupled dynamic units whose states get excited upon a particular local event. It has the capacity to propagate waves of excitation, which often self-organize into spiral patterns. Examples of such systems in nature include neural networks, Belousov-Zhabotinsky reaction,  as well as coupled oscillators such as fireflies and pacemaker cells. In a discrete setting, excitable media can be modeled using the framework of generalized cellular automaton (GCA). Given a simple connected graph $G=(V,E)$ and a fixed integer $\kappa\ge 2$,  the microstate of the system at a given discrete time $t\ge 0$ is given by a $\kappa$-coloring of vertices $X_{t}:V\rightarrow \mathbb{Z}_{\kappa}=\mathbb{Z}/\kappa\mathbb{Z}$. A given initial coloring $X_{0}$ evolves in discrete time via iterating a fixed deterministic transition map $\tau: X_{t}\mapsto X_{t+1}$, which depends only on local information at each time step. This generates a trajectory $(X_{t})_{t\ge 0}$, and its  limiting behavior in relation to the topology of $G$ and structure of $\tau$ is of our interest. 

\qquad Greenberg-Hastings Model (GHM) and cyclic cellular automaton (CCA) are two particular GCA models for excitable media  \cite{wiener1946mathematical} that have been studied extensively in the 1990's. GHM was introduced by Greenberg and Hastings \cite{greenberg1978spatial} to capture phenomenological essence of neural networks in a discrete setting, whereas CCA was introduced by Bramson and Griffeath \cite{bramson1989flux} as a discrete time analogue of the cyclic particle systems. In GHM, think of each vertex of a given graph as a $\kappa$-state neuron. An excited neuron (i.e., one in state 1) excites neighboring neurons at rest (i.e., in state 0) and then needs to wait for a refractory period of time (modeled by the remaining $\kappa-2$ states) to become rested again. In CCA, each vertex of the graph is inhabited by one of $\kappa$ different species in a cyclic food chain. Species of color i are eaten (and thus replaced) by species of color $(i+1)$ mod $\kappa$ in the neighborhood at each time step. More precisely, the transition maps  $X_{t}\mapsto X_{t+1}$ for $\kappa$-color GHM and CCA are given below: 
\begin{equation}
	(\text{GHM})\qquad X_{t+1}(v)=\begin{cases}
1 & \text{if $X_{t}(v)=0$ and $\exists u\in N(v)$ s.t. $X_{t}(u)=1$  }\\
0 & \text{if $X_{t}(v)=0$ and $\nexists u\in N(v)$ s.t. $X_{t}(u)= 1$  }\\
X_{t}(v)+1\,\,\text{(mod $\kappa$)} & \text{otherwise}
\end{cases}
\end{equation}
\begin{equation}
	(\text{CCA})\qquad X_{t+1}(v)=\begin{cases}
X_{t}(v)+1 \,\, (\text{mod $\kappa$})& \text{if $\exists u\in N(v)$ s.t. $X_{t}(u)=X_{t}(v)+1$ (mod $\kappa$) }\\
X_{t}(v) & \text{otherwise}
\end{cases}
\end{equation}

\qquad GHM and CCA are among the few models of discrete excitable media that have been studied rigorously in the probability literature, mostly on the integer lattice $\mathbb{Z}^{d}$ with randomly chosen initial $\kappa$-coloring $X_{0}$. For $d=1$, Fisch \cite{fisch1990one} showed that CCA exhibits a phase transition ``between'' $\kappa=4$ and 5; for $\kappa\in\{3,4\}$ each vertex increments its color infinitely often with probability 1, whereas for $\kappa\ge 5$ the dynamics fixates. In a subsequent work \cite{fisch1992clustering}, Fisch established clustering of the one-dimensional 3-color CCA by showing that the density of borders between different colors approaches 0 at the rate $t^{-1/2}$ as time $t$ increases.  The main technique was a connection to a random walk, which was adapted to the 3-color GHM on $\mathbb{Z}$ by Durrett and Steif \cite{durrett1991some}, and to GHM on $\mathbb{Z}$ with arbitrary $\kappa$ by Fisch and Gravner \cite{fisch1995one}. 

\qquad In higher dimensions, waves of excitation can propagate feed back on itself, resulting in self-sustained local wave generators, the  stable periodic objects (SPOs). Fisch, Gravner, and Griffeath showed that the limiting behavior of CCA on $\mathbb{Z}^{d}$ for any $d\ge 2$ and $\kappa\ge 3$ is governed by the formation of SPOs  \cite{fisch1991cyclic}. The 3-color GHM on $\mathbb{Z}^{d}$ shares this behavior \cite{durrett1991some}, and similar behavior of GHM and CCA for any $\kappa\ge 3$ on higher dimensions was studied by Fisch and Gravner \cite{fisch1991threshold}, who introduced an additional parameter $\theta$, the threshold number of excited neighbors required to excite a vertex. 

\qquad A fundamental difficulty in understanding the limiting behavior of excitable media models on general graphs is the complexity in generation of SPOs and interactions between them. On trees, however, these objects are topologically prohibited. This makes trees special, and a substantial portion of the paper addresses the behavior of CCA and GHM on them. Related models on trees arose in the context of distributed algorithms for digital clock synchronization \cite{dolev2000self}: a 3-color GCA model on finite trees studied by Herman and Ghosh \cite{herman1995stabilizing}; odd $\kappa\ge 3$ models, which coincide with the 3-color CCA in case $\kappa=3$, investigated by Boulinier, Petit, and Villain \cite{boulinier2006toward}; and a GCA model for pulse-coupled inhibitory oscillators, called the firefly cellular automaton, whose behavior on finite trees was recently addressed by the second author \cite{lyu2015synchronization}.

\qquad In this paper, we revisit the 3-color GHM and CCA and characterize their limiting behavior on arbitrary graphs. Our main technique is the construction of a monotone comparison process on the universal covering space of the underlying graph. This comparison process itself is inspired by a famous consensus algorithm of Lamport \cite{lamport1978time}: if the vertices of a connected graph are equipped with integer-valued opinions (not necessarily distinct), and at each step each node simultaneously adopts the maximum opinion among itself and its neighbors. Then in some finite time the entire graph reaches a consensus, which is the initial global maximum. The correspondence between the GHM and CCA dynamics and this comparison process is such that a vertex is excited in the original dynamics if and only if all vertices in its fiber increment their opinion by 1. We remark that a similar observation was made by Belitsky and Ferrari \cite{belitsky1995ballistic} for one-dimensional ballistic annihilation system. A consequence of this comparison for the 3-color GHM and CCA on arbitrary graphs is that the SPOs are static and cannot arise spontaneously. After establishing the comparison process, we apply it to deterministic and random finite graphs and to infinite trees.

\section{Statements of results}\label{section:results}

Let $G=(V,E)$ be a connected graph and let $X_{0}:V\rightarrow \mathbb{Z}_{3}$ be an initial 3-coloring, and let $(X_{t})_{t\ge 0}$ be the resulting CCA or GHM dynamics. We say a node $x$ is \textit{excited} at time $t$ if: $X_{t+1}(x) = X_t(x)+1\mod 3$ in the case of CCA; and if $X_{t}(x) = 0$ and   $X_{t+1}(x)= 1$ in the case of GHM. Define $\nbe_{t}(x)$ to be the number of excitations $x\in V$ undergoes in the first $t$ iterations,
	\begin{equation}
	\nbe_{t}(x) = \sum_{s=0}^{t-1} \indicator(\text{$x$ is excited at time $s$}),\qquad \text{$x\in V$ and $t\ge 1$,}  
	\end{equation}
where $\indicator(A)$ denotes the indicator function of event $A$.

\qquad Let $\overline{E}$ be set of all ordered pairs of adjacent nodes, i.e., $\overline{E} = \{(u,v)\in V^{2}\,|\, \text{$uv\in E$}\}$. Let $X:V\rightarrow \mathbb{Z}_{3}$ be a 3-coloring on $G$. For the CCA, we define the associated anti-symmetric 1-form $dX:\overline{E}\rightarrow \{-1,0,1\} $ as follows:
\begin{equation}
	dX(u,v) = X(v)-X(u)
\end{equation}

where the subtraction is taken in $\mathbb{Z}_{3}$. Similarly, for the GHM dynamics, we define 
\begin{equation}
	dX(u,v) =
	\begin{cases}
		1 & \text{if $X(u)=0$ and $X(v)=1$}\\
		-1 & \text{if $X(u)=1$ and $X(v)=0$}\\
		0 & \text{otherwise.}
	\end{cases}	
\end{equation}
The definition is such that a vertex $x\in V$ is excited at time $t$ if and only if $dX_{t}(x,y)=1$ for some $y\in N(x)$, in which case we say that $y$ {\it excites\/} $x$ at time $t$.

\qquad A \textit{walk} $\vec{W}$ in $G$ is a finite or infinite sequence of nodes $(v_{i})_{i\ge 0}$ in $G$ such that $v_{i}$ is adjacent to $v_{i+1}$ for all $i\ge 0$.  We say that $\vec{W}$ is \textit{non-backtracking} if $v_{i}\ne v_{i+2}$ for all $i\ge 0$, and \textit{closed} if it consists of finitely many nodes where the first and last nodes coincide, and a \textit{directed path} if all nodes in $\vec{W}$ are distinct. We say $\vec{W}$ is a \textit{cycle}  if it is a closed walk in which only first and last vertex agree.  For a walk $\vec{W}=(v_{i})_{i\ge 0}$ and a 3-coloring $X:V\rightarrow \mathbb{Z}_{3}$, we define the \textit{path integral} of $X$ on $\vec{W}$ by 
\begin{equation}
	\int_{\vec{W}} \,dX := \sum_{i\ge 0} dX(v_{i},v_{i+1}) 
\end{equation}

We say $dX$ (with respect to CCA or GHM dynamics) is \textit{irrotational} if every contour integral of $X$ over directed cycles is zero, that is, 
	\begin{equation}\label{cauchyint}
		\oint_{\vec{C}} \,dX \equiv 0 	
	\end{equation}
for all closed directed cycles $\vec{C}$ in $G$. Note that any closed walk can be decomposed into a finite number of directed cycles. Hence if $dX$ is irrotational, then (\ref{cauchyint}) holds for all closed walks~$\vec{W}$.

\qquad One of our key tools is the following lemma:

\begin{customlemma}{1}\label{key_lemma}
	Let $G=(V,E)$ be an arbitrary connected graph and let $(X_{t})_{t\ge 0}$ be a CCA or GHM trajectory of an initial 3-coloring $X_{0}:V\rightarrow \mathbb{Z}_{3}$. Then for each $x\in V$ we have 
	\begin{equation}
		\nbe_{t}(x) = \max_{|E(\vec{W})|\le t} \int_{\vec{W}} \,dX_{0},
	\end{equation}
where the maximum runs over all walks $\vec{W}$ starting from $x$ of length $|E(\vec{W})|$ at most $t$. 
\end{customlemma}

For each $x\in V$, we define the following quantity 
\begin{equation}
\alpha(x) = \limsup_{t\rightarrow \infty} \frac{\nbe_{t}(x)}{t}.
\end{equation}

An immediate consequence from the above lemma is the following:

\begin{customcorollary}{2}\label{activity}
	Let $G=(V,E)$ and $(X_{t})_{t\ge 0}$ be as before. If $\sup_t\nbe_{t}(x)<\infty$ holds for some $x\in V$, then it holds 
	for all $x\in V$.	Moreover, for each $x,y\in V$, we have $\alpha(x)=\alpha(y)$. 
\end{customcorollary}

\qquad We say that $X_{t}$ \textit{fixates} if $\nbe_{t}(x)$ is bounded in time for some (an thus for all) 
$x\in V$, and that it \textit{fluctuates} otherwise. We say $X_{t}$ \textit{synchronizes} if for every two vertices $x,y\in V$, there exists $N=N(x,y)\in \mathbb{N}$ such that $X_{t}(x)=X_{t}(y)$ for all $t\ge N$.  It is not hard to see that fixation and synchronization are equivalent notions for 3-color CCA and GHM dynamics. In fact, this is true for the general $\kappa$-color GHM dynamics, but not for $\kappa$-color CCA dynamics with $\kappa\ge 4$ as there are many non-interacting pairs of colors. 
Furthermore, from now on we will denote the constant value of $\alpha(x)$ by $\alpha$, and call it 
the \textit{activity} of the dynamics $(X_{t})_{t\ge 0}$.  Finally, we say that $X_{t}$ \textit{synchronizes weakly} if $\alpha=0$ and that it \textit{oscillates} otherwise.

\qquad Our first result is a characterization of the limiting behavior of the 3-color CCA or GHM dynamics on finite graphs in terms of the irrotationality of the induced 1-form at time 0:

\begin{customthm}{3}\label{mainthm1}
	Let $G=(V,E)$ be a finite connected graph and let $(X_{t})_{t\ge 0}$ be a CCA or GHM trajectory of an initial 3-coloring $X_{0}:V\rightarrow \mathbb{Z}_{3}$. Then $X_{t}$ synchronizes if and only if $dX_{0}$ is irrotational. Furthermore, we have the following:
	
\begin{description}[noitemsep]
	\item{(i)} If $dX_{0}$ is irrotational, then $X_{t}(x)=X_{t}(y)$ for all $x,y\in V$ and $t\ge D$, where $D$ is the diameter of $G$;
	\item{(ii)} If $dX_{0}$ is not irrotational, then for each node $x\in V$, we have 
	\begin{equation}\label{finite_activity}
		 \alpha=\lim_{t\rightarrow \infty} \frac{\nbe_{t}(x)}{t} =  \max_{\vec{C}} \frac{1}{|V(\vec{C})|} \oint_{\vec{C}} dX_{0}  
	\end{equation}
where the maximum runs over all closed directed cycles $\vec{C}$ in $G$. 
\end{description}

\end{customthm}  

One immediately has the following corollary:

\begin{customcorollary}{4}\label{corollary1}
	Let $G=(V,E)$ be a finite graph. Then under either CCA or GHM dynamics, arbitrary 3-coloring $X_{0}:V\rightarrow \mathbb{Z}_{3}$ on $G$ synchronizes if and only if $G$ is a tree.  
\end{customcorollary}

\qquad We next present a result on the CCA and GHM dynamics on general graphs $G=(V,E)$, starting from a random 3-coloring $X_{0}$ drawn from the uniform product measure on $\mathbb{Z}_{3}^{V}$, which generates a random 3-color CCA or GHM trajectory. A \textit{matching} in $G$ is a set of edges $\{e_{1},\cdots,e_{k}\}$ $\subset E$ where two edges do not share a common vertex.

\begin{customthm}{5}\label{mainthm2_cycle}
	Consider a random 3-color CCA or GHM trajectory $(X_{t})_{t\ge 0}$ on a graph $G=(V,E)$. If $G$ contains a cycle, $X_t$ oscillates with a positive probability. Furthermore, suppose $G$ has a matching $\{e_{1},\cdots,e_{k}\}$ and distinct cycles $C_{1},\cdots,C_{k}$ (not necessarily vertex-disjoint) such that $e_{i}\in E(C_{j})$ iff $i=j$ for all $1\le i\le j\le k$. Then 
	\begin{equation}
		\mathbb{P}(\text{$X_{t}$ synchronizes weakly})\le (7/9)^{k} .
	\end{equation}
This bound is achieved when $C_{i}$'s are vertex-distjoint triangles and assuming CCA dynamics. In particular, if $k=\infty$, i.e., there is an infinite matching and an infinite sequence of cycles as above, then $X_{t}$ oscillates almost surely.
\end{customthm}

An interesting application of the above theorem is the following phase transition of limiting behavior of $X_{t}$ on the Erd\"os-R\'enyi random graph.

\begin{customthm}{6}\label{corollary_ER}
	Let $G=G(n,p)$ be the Erd\"os-R\'enyi random graph and let $(X_{t})_{t\ge 0}$ be a random CCA or GHM trajectory, and let $\mathbb{P}$ denote the joint product probability measure for $(G,X_0)$.
	\begin{description}[noitemsep]
		\item{(i)} If $p=o(1/n)$ then $X_{t}$ synchronizes on each component of $G$ a.a.s.;
		\item{(ii)} If $p=\lambda/n$ for any $0<\lambda<1$, then there exists a constant $C=C(\lambda)\in (0,1)$ such that 
\begin{equation}
	\lim_{n\rightarrow \infty}\mathbb{P}(\text{$X_{t}$ synchronizes on each component of $G(n,p)$})  =  C.
\end{equation}				
 		\item{(iii)} If $p=\lambda/n$ for any $\lambda>1$, then there exists a constant $D=D(\lambda)>0$ such that for all sufficiently large $n$, 
 		\begin{equation}
 			\mathbb{P}(\text{$X_{t}$ oscillates on the largest component of $G(n,p)$})\ge 1-e^{-Dn}.
 		\end{equation}
 	\end{description}
\end{customthm}
Explicit expressions for the constants $C$ and $D$ in the assertion are given in the proof.

 \begin{figure*}[h]
 	\centering
 	\includegraphics[width=0.92 \linewidth]{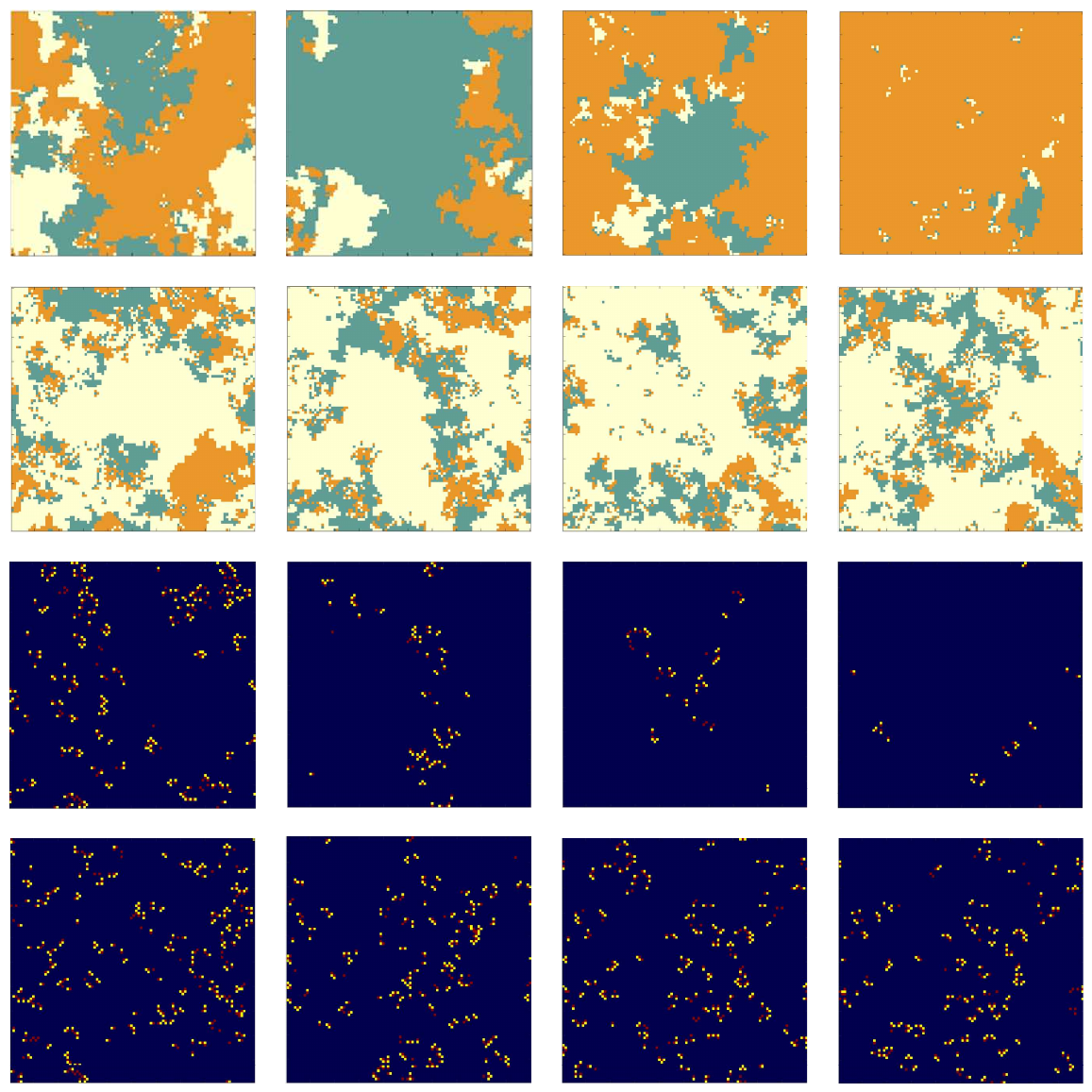}
 	\caption{ (Top row) Snapshots of 3-color CCA on a uniform spanning tree of a 100 by 100 torus, each 100 iterations from left to right. (Second row) Dynamics after 12 random edges are added to the spanning tree. Orange =0, green=1, and yellow=2. Corresponding simulations for 3-color GHM are shown in the third and fourth rows. Dark blue=0, yellow=1, and red=2.
 	}
 	\label{snapshots}
 \end{figure*}

 \qquad We next address the limiting behavior on infinite trees. Let $\Gamma=(V,E)$ be an infinite, but locally finite, tree rooted at $0\in V$. In the context of trees, $X_{0}$ will be a random 3-coloring on the vertices drawn from a general product measure $\mathbb{P}=\mathbb{P}_{p_{0},p_{1},p_{2}}$ with marginal distribution $\mathbb{P}(X_{0}(x)=i)=p_{i}$ for all $x\in V$ and $i\in \mathbb{Z}_{3}$.  The reason we are considering this general situation is to highlight an intriguing symmetry between the three colors on infinite trees; by contrast, in Theorems~\ref{mainthm2_cycle} and~\ref{corollary_ER} the added generality provides little of interest. The random coloring $X_{0}$ induces random variables $dX_{0}(x,y)$ on each adjacent pair $(x,y)\in \overline{E}$, which we will refer to as {\it increments\/}. Note that the increments are bounded, have zero expectation, and are identically distributed. Also note that they are in general not independent. For instance, no consecutive increments can be both 1 in case of GHM, for arbitrary marginal density. However, it is easy to see that in our case the increments are \textit{1-correlated}, by which we mean that two increments are independent if their underlying edges are vertex disjoint.

\qquad We append a random variable $S_{\sigma}$ at each vertex $\sigma\in V$ by 
\begin{equation}
	S_{\sigma} = \int_{\vec{P}_{\sigma}} dX_{0}
\end{equation} 
where $\vec{P}_{\sigma}$ is the unique directed path from $0$ to $\sigma$. The collection $\{S_{\sigma}\}_{\sigma\in V}$ is called a \textit{$\Gamma$-indexed walk}. If we take $\Gamma$ to be a Galton-Watson tree, for instance, then a $\Gamma$-indexed walk can be viewed as a branching random walk. In the beginning of Section 5 we will see that the activity $\alpha$ of the dynamics $(X_{t})_{t\ge 0}$ coincides with the following quantity 
\begin{equation}
	v_{c}:=\limsup_{n\rightarrow \infty}\frac{1}{n} \max_{|\sigma|=n}S_{\sigma}, 
\end{equation}     
which is called the \textit{cloud speed} of the associate $\Gamma$-indexed walk. Hence on infinite trees,  understanding the CCA and GHM dynamics boils down to the study of the cloud speed of tree indexed walks. Hence we will be interested in the cloud speed of a $\Gamma$-indexed random walk, where $\Gamma$ is arbitrary and increments are 1-correlated.

\qquad The value of $v_{c}$ is closely related to the exponential growth rate of populations in $\Gamma$. Perhaps the simplest quantity to measure the growth is its
\textit{volume entropy} given by 
\begin{equation}
	\entropy(\Gamma)=\limsup_{n\rightarrow \infty}\frac{1}{n} \log A_{n}, 
\end{equation}
where $A_{k}$ is the number of vertices of $\Gamma$ at level $k$. We remark that replacing $A_{k}$ by the number $B_{k}$ of vertices of $\Gamma$ upto level $k$ in the above definition does not change the value, which justifies its name. (This quantity is the \textit{Minkowski dimension} of the boundary space of $\Gamma$; see e.g., \cite{benjamini1994tree}, which considers only trees without leaves.) Consider first $\Gamma=\mathbb{Z}$, in which case $\entropy(\Gamma)=0$. The law of iterated logarithm says that for random walks with centered i.i.d. increments, the maximum grows in the order of $\sqrt{n\log\log n}$ a.s., so $v_{c}=0$ a.s. More generally, positive volume entropy is necessary for positive cloud speed; see Theorem \ref{mainthm2} (iii).

\qquad However, $\entropy(\Gamma)>0$ is not sufficient to guarantee $v_{c}>0$, which we illustrate in the following example. Let $\Gamma^{a}$ be obtained from a single infinite ray $\gamma$ by attaching $2^{|\sigma|}-1$ leaves to each vertex $\sigma\in V(\gamma)$. Even though the $k^\text{th}$ level contains $2^{k}$ vertices so that $\entropy(\Gamma)=\log 2$, on this tree the cloud speed is zero, as there is no contribution of the increments on leaves to the cloud speed. Thus the volume entropy may not be a suitable notion to study the cloud speed. Nonetheless, under the assumption that $\Gamma$ has no leaves and the increments of the $\Gamma$-indexed walk are centered, i.i.d., and satisfy a mild moment condition, Benjamini and Peres \cite{benjamini1994tree} showed that $v_{c}>0$ if and only if $\entropy(\Gamma)>0$. Furthermore, they obtained sharp upper and lower bounds on $v_{c}$ in terms of the volume entropy and the large deviations rate of one-dimensional random walk with the same step distribution. 

\qquad The example $\Gamma^{a}$ illustrates that for trees with leaves we may need a refined quantity that measures the average number of infinite branches per vertex. Lyons \cite{lyons1990random} introduced such a quantity called the \textit{branching number}. Given a tree $\Gamma$, a \textit{cutset} $\Pi$ is a finite set of vertices not including 0 such that every infinite path from $0$ intersects $\Pi$ and such that there is no pair $\sigma,\tau\in \Pi$ with $\sigma<\tau$. The \textit{branching number} of $\Gamma$ is defined by 
\begin{equation}
	\br(\Gamma) = \inf\left\{  \lambda>0 \,\bigg|\, \inf_{\Pi} \sum_{\sigma\in \Pi} \lambda^{-|\sigma|}=0 \right\}. 
\end{equation}

It follows from the definitions that $\log \br(\Gamma)\le \entropy(\Gamma)$, and it is known that equality holds when $\Gamma$ is sufficiently regular, e.g., for almost all trees generated by a Galton-Watson process \cite{lyons1990random}. Turning back to the example $\Gamma^{a}$, indeed we have $\log \br(\Gamma^{a})=0$ since there is only one infinite branch, which seems to correspond with the fact that $v_{c}=0$ on $\Gamma^{a}$. However, it turns out that the branching number does not decide when the cloud speed vanishes either, as there are trees with branching number 1 but large cloud speed. We will give such an example in Section 5 (Example \ref{geometric tree})

\qquad This leads us to introduce another notion of dimension of $\Gamma$, with a better connection to the cloud speed in case $\Gamma$ has leaves, which allows us to improve known lower and upper bounds. Given a tree $\Gamma$ and for each $n\le m$, denote by $A_{n,m}$ the number of vertices at level $n$ which have descendants at level $m$. For each $r> 1$, define the \textit{r-volume entropy} of $\Gamma$ by 
\begin{equation}\label{r_entropy_def}
	\entropy_{r}(\Gamma) =  \limsup_{n\rightarrow  \infty}\frac{1}{n} \log A_{n,\floor{rn}},
\end{equation}
and its \textit{reduced volume entropy} by 
\begin{equation}
	\rentropy(\Gamma) =  \limsup_{r\searrow 1} \entropy_{r}(\Gamma).
\end{equation}

In particular, if $\Gamma$ has no leaves, then $A_{n}=A_{n,\floor{rn}}$ for all $r> 1$ so $\entropy_{r}(\Gamma)=\entropy(\Gamma)$ for all $r> 1$. In general, we have $\log \br(\Gamma)\le  \entropy_{r}(\Gamma)\le \rentropy(\Gamma)\le \entropy(\Gamma)$ for all $r>1$.

\qquad In the statement of our main theorem to follow, recall that we consider a product measure with unequal densities $p_0$, $p_1$, $p_2$ of colors $0$, $1$, $2$. In preparation, we introduce certain functions that determine the relevant large deviation rates.  
For each $t\ge 0$, let $x=x(t)$ be the largest positive root
of the cubic equation 
\begin{equation}\label{defining_cubic_CCA}
x^{3}-x^{2}=p_{0}p_{1}p_{2}(e^{3t}+e^{-3t}-2),
\end{equation}
and define $\Lambda_{\text{CCA}}(t)=\log x(t)$.
Note that $x(t)$ is strictly increasing on $[0,\infty)$, as the right-hand side of (\ref{defining_cubic_CCA}) is nonnegative and strictly increasing for $t\in [0,\infty)$, and the 
left-hand side strictly increases from $0$ to $\infty$ for $x\in[1,\infty]$. 
Hence 
$\Lambda_{\text{CCA}}(t)\ge 0$ is well-defined for $t\in [0,\infty)$, strictly increasing and continuous.   Next, let $\Lambda_{\text{CCA}}^{*}(u)=\sup_{t\ge 0} \{ut-\Lambda_{\text{CCA}}(t)\}$ be its Legendre transform. Observe that $x(t)\sim (p_0p_1p_2)^{1/3}e^{t}$ as $t\to\infty$, so that 
\begin{equation}
\Lambda_{\text{CCA}}(t)=t+\frac{1}{3}\log (p_{0}p_{1}p_{2}) +o(1). 
\end{equation}
Hence it follows that $\Lambda_{\text{CCA}}^{*}(u)$ is finite, strictly increasing and continuous 
on $[0,1]$, $\Lambda_{\text{CCA}}^{*}(1)=-\frac{1}{3}\log (p_{0}p_{1}p_{2})$, and  
$\Lambda_{\text{CCA}}^{*}\equiv\infty$ on $(1,\infty)$. Moreover,
it is easy to check that for $p_{0}=p_{1}=p_{2}=1/3$, we have the following closed-form expressions: 
\begin{eqnarray}
&&\Lambda_{\text{CCA}}(t)=\log(1+e^{t}+e^{-t})-\log 3\\
&&\Lambda_{\text{CCA}}^{*}(u) = 
\begin{cases}
u\log\left( \frac{u+\sqrt{4-3u^{2}}}{2(1-u)} \right) - \log \left( \frac{1+\sqrt{4-3u^{2}}}{1-u^{2}}\right)+\log 3 & \text{if $u\in [0,1)$}\\
\log 3 & \text{if $u=1$}  \\
\infty & \text{if $u\in (1,\infty)$}
\end{cases}
\end{eqnarray}  
Finally, define 
$\Lambda_{\text{GHM}}(t)=\Lambda_{\text{CCA}}(t/3)$ and let $\Lambda_{\text{GHM}}^{*}$ be its Legendre transform, so that $\Lambda_{\text{GHM}}^{*}(u)=\Lambda_{\text{CCA}}(3u)$. We now 
state our main theorem.

\begin{customthm}{7}\label{mainthm2}
Let $\Gamma=(V,E)$ be an infinite rooted tree and $(X_{t})_{t\ge 0}$ the random 3-color CCA or GHM trajectory on $\Gamma$, where $X_{0}$ is drawn from the product measure with marginal density $\mathbb{P}(X_{0}(\sigma)=i)=p_{i}$ for $i\in \mathbb{Z}_{3}$. Denote by $\Lambda$ either $\Lambda_{\text{CCA}}$ or 
$\Lambda_{\text{GHM}}$, depending on the dynamics, and let $\Lambda^{*}$ be corresponding Legendre transform.  Then we have the following:
\begin{description}[noitemsep]
	\item{(i)} $X_{t}$ synchronizes weakly if and only if $\,\rentropy(\Gamma)=0$.   
	\item{(ii)} The activity $\alpha$ is equal to the cloud speed $v_{c}$ of the associated $\Gamma$-indexed random walk $\{S_{\sigma}\}_{\sigma\in V}$.	
	\item{(iii)} $v_{c}$ satisfies the  upper bound 
	\begin{equation}\label{vc_upper_entropy_thm2}
		\Lambda^{*}(v_{c})\le \rentropy(\Gamma).
	\end{equation}
		\item{(iv)}  Suppose $\log \br(\Gamma)=\entropy(\Gamma)$ and let $B=1$ for CCA and $B=1/3$ for GHM. If $\br(\Gamma)\ge 1/\sqrt[3]{p_{0}p_{1}p_{2}}$, then $v_{c}=B$.  Otherwise, $v_{c}<B$ and  equality holds in (\ref{vc_upper_entropy_thm2}), which determines $v_c$. In particular, $v_{c}$ for GHM is a third of that for CCA, for arbitrary $p_{0}$, $p_1$, $p_2$.
	\item{(v)}  For each $r>1$ such that $v_c/(1-r^{-1}) < B$, where $B$ is defined in (iv), we have 
	\begin{equation}\label{vc_lowerbd_rcut}
		\frac{\entropy_{r}(\Gamma)}{r-1}\le \Lambda^{*}\left( \frac{v_{c}}{1-r^{-1}}\right).
\end{equation}		
In particular, if $\Gamma$ has no leaves, then 
	\begin{equation}
		 \Lambda\left(h(\Gamma)/v_{c}\right)\le \entropy(\Gamma),
	\end{equation}	
	which is sharp in the sense that there exists a tree-indexed walk on some tree $\Gamma^{-}$ with cloud speed $v_{c}^{-}$ satisfying $\Lambda(\entropy(\Gamma)/(v_{c}^{-}-\epsilon))>\entropy(\Gamma)$ for each $\epsilon>0$.\\
\end{description}  
\end{customthm}

\begin{customremark}{8}	
	Note that Theorem \ref{mainthm2} (iv) implies that, when $\log \br(\Gamma)=\entropy(\Gamma)$, which for instance holds a.s. for Galton-Watson trees, the activity $\alpha$ for both CCA and GHM depends only on the product $p_{0}p_{1}p_{2}$; in 
	particular, it is symmetric in the three densities. This is particularly surprising for GHM, in which the three colors play very different roles. 
	For example, pick a small $\epsilon>0$ and compare two initial states: one with $p_0=p_1=(1-\epsilon)/2$ and $p_2=\epsilon$, and the other with $p_0=p_2=(1-\epsilon)/2$ and $p_1=\epsilon$. The two generate very different dynamics; nevertheless, the higher annihilation rate of the first one causes their activities to match. 
\end{customremark}

\qquad This paper is organized as follows. In Section 3 we prove Lemma \ref{key_lemma} by introducing a monotone comparison process defined on the universal covering space of the underlying graph. In Section 4 we give proofs for the main theorems concerning finite and random graphs, Theorem \ref{mainthm1}, \ref{mainthm2_cycle}, and \ref{corollary_ER}. In Section 5, we extend Proposition 4.1 in Benjamini and Peres \cite{benjamini1994tree} and give sharp estimates of the cloud speed of a tree-index random walk when the underlying tree is arbitrary and the increments are 1-correlated. By the method developed in Section 3, this leads to the proof of Theorem \ref{mainthm2}.

\section{Tournament expansion of 3-color GHM and CCA}

We prove Lemma \ref{key_lemma} in this section.  One of the main complication in understanding the 3-color CCA and GHM dynamics is that the color space $\mathbb{Z}_{3}$  has a cyclic hierarchy so that the configuration space lacks monotonicity. To overcome this difficulty, we introduce a simple monotone comparison process and establish its relationship with our CCA and GHM dynamics.

\qquad \textit{Tournament process} on a graph is a simple deterministic process in which initially all nodes have an integer rank, and in each step each node simultaneously adopts the maximum rank among itself and its neighbors.  More precisely, given a connected graph $G=(V,E)$, a map $\mathtt{rk}_{0}:V\rightarrow \mathbb{Z}$ is called a \textit{ranking} on $G$. The transition map from time $t$ to $t+1$ is given by 
\begin{equation}
	\mathtt{rk}_{t+1}(x) = \max \{\mathtt{rk}_{t}(y)\,|\, y\in N(x)\cup\{x\}\}. 
\end{equation} 
Iteration of the above transition rule generates a discrete-time orbit $(\mathtt{rk}_{t})_{t\ge 0}$ of rankings.

\qquad Observe that if $G$ is finite,  then for any initial ranking $\mathtt{rk}_{0}$ on G there is a global maximum, which every node will eventually achieve. In general, locally maximum rank propagates with unit speed across the graph until it is overcome by waves from higher rankers. To make this observation precise, or each node $x\in V$ and radius $t\in \mathbb{N}$, define $M_{t}(x)$ to be the maximum initial rank in the $t$-ball centered at $x$:
\begin{equation}
	M_{t}(x) = \max\{\mathtt{rk}_{0}(y)\,|\, d(x,y)\le t\} 
\end{equation} 
where $d(x,y)$ is the usual graph distance in $G$. The following proposition characterizes the dynamics in tournament processes:

\begin{customprop}{3.1}\label{tournamentdynamics}
Let $G=(V,E)$ be a connected graph and fix an initial ranking $\mathtt{rk}_{0}:V\rightarrow \mathbb{Z}$. Then the local dynamic at any given node $x\in V$ is given by 
\begin{equation}
	\mathtt{rk}_{t}(x) = M_{t}(x)
\end{equation}
\end{customprop}

\begin{proof}
	Fix $x\in V$ and $t\ge 1$. By tracing back the origins of ranks, it is easy to see that if $\mathtt{rk}_{t}(x)=r>\mathtt{rk}_{0}(x)$, then necessarily there exists $y\in V$ with $d(x,y)\le t$ and $\mathtt{rk}_{0}(y)=r$. This yields $\mathtt{rk}_{t}(x)\le M_{t}(x)$. On the other hand, let $y\in V$ be the vertex that attains $M_{t}(x)$, i.e., $d(x,y)\le t$ and $\mathtt{rk}_{0}(y)=M_{t}(x)$. Choose a shortest path $P$ from $x$ to $y$. Since $P$ has length $\le t$, the rank $\mathtt{rk}_{0}(y)$ of $y$ propagates along $P$ and reaches $x$ by time $t$. Since $\mathtt{rk}_{s}(x)$ is non-decreasing in $s$, this yields that $\mathtt{rk}_{t}(x)\ge M_{t}$. This shows the assertion. 
\end{proof}

\qquad Having introduced the tournament process, we now establish a comparison tournament process for CCA or GHM dynamics. Fix a graph $G=(V,E)$ and an initial 3-coloring $X_{0}:V\rightarrow \mathbb{Z}_{3}$. Designate an arbitrary vertex $x$ in $G$ as its base point. If $xy\in E$ and $dX_{0}(x,y)=1$, this means $y$ excites $x$ at time $0$ and $X_{1}(x)=X_{0}(y)$. Hence when we construct an associated tournament process, it would be natural to give $y$ a rank that is greater than that of $x$, e.g., $\mathtt{rk}_{0}(y)=\mathtt{rk}_{0}(x)+dX_{0}(x,y)$. In general, if $z$ is any node in $G$, we may try to define the rank of $z$ by 
\begin{equation}\label{initialranking_def}
	\mathtt{rk}_{0}(z)=\mathtt{rk}_{0}(x)+\int_{\vec{P}} dX_{0}
\end{equation}
where $\vec{P}$ is any walk from $x$ to $z$ with initial condition $\mathtt{rk}_{0}(x)=0$. However, this gives a well-defined ranking on $G$ if and only if $dX_{0}$ is irrotational. In particular, $\mathtt{rk}_{0}$ is well-defined when $G$ is a tree. In that case, we define the \textit{tournament expansion} of $(X_{t})_{t\ge 0}$ to be the tournament process $(\mathtt{rk}_{t})_{t\ge 0}$.  

\qquad To make sense of the above definition on general graphs, we need to distinguish the endpoints of distinct walks from $x$. This encourages us to define the associated tournament process on the universal covering space of $G$. The \textit{universal covering space of $G$ with base point $x$} is a tree $\mathcal{T}_{x}=(\mathcal{V},\mathcal{E})$ where $\mathcal{V}$ is the set of all finite non-backtracking walks in $G$ starting from $x$, and adjacency is given by one-step extension: $\{\vec{W},\vec{W}'\} \in \mathcal{E}$ if and only if one of them can be obtained by adjoining a single vertex at the end of the other walk. We have the natural covering map $p:\mathcal{V}\rightarrow V$ which maps each finite  walk $\vec{W}$ to its endpoint. We identify the length 0 walk from $x$ to $x$ with $x\in V$. It is customary to denote $p:\mathcal{T}_{x}\rightarrow G$. For each $y\in V$, $\tilde{y}$ denotes any node in $\mathcal{T}_{x}$ such that $p(\tilde{y})=y$. For each 3-coloring $X$ on $G$, we denote by $\tilde{X}$ the \textit{lift} of $X$ onto $\mathcal{T}_{x}$, defined by $\tilde{X}(\tilde{x})=X(x)$ for all $\tilde{x}\in p^{-1}(x)$. A trivial but important observation is that any GCA dynamics of $(G,X_{0})$ naturally lifts to $(\mathcal{T}_{x},\tilde{X}_{0})$, due to the fact that the transition map is completely determined locally. Hence, in particular, we can lift any CCA or GHM dynamics on any simple graph $G$ onto its universal covering space. 

\qquad Now for a given 3-color CCA or GHM dynamics $(X_{t})_{t\ge 0}$ on arbitrary connected graph $G=(V,E)$ with base point at $x\in V$, we define its \textit{tournament expansion} by the tournament expansion of $(\tilde{X}_{t})_{t\ge 0}$, the 3-color CCA or GHM dynamics lifted onto the universal covering space $\mathcal{T}_{x}$ with base point at $x\in \mathcal{V}$. The following lemma establishes that the two dynamical systems are compatible: 

\begin{customlemma}{3.2}\label{tournament_lemma}
	Let $G=(V,E)$, $\mathcal{T}_{x}=(\mathcal{V},\mathcal{E})$, $(X_{t})_{t\ge 0}$, and $(\mathtt{rk}_{t})_{t\ge 0}$ be as before. Then for all $t\ge 0$ and $\tilde{z}\in \mathcal{V}$, we have   
	\begin{equation}
	\mathtt{rk}_{t+1}(\tilde{z})-\mathtt{rk}_{t}(\tilde{z}) = \indicator(\text{$z$ is excited at time $t$}).
	\end{equation}
	In particular, for all $t\ge 0$, we have 
	\begin{equation}
		\mathtt{rk}_{t}(x) = \nbe_{t}(x).
	\end{equation}
\end{customlemma} 
Note that the key lemma (Lemma \ref{key_lemma}) follows from Lemma \ref{tournament_lemma} in conjunction with Proposition \ref{tournamentdynamics}. We prove Lemma \ref{tournament_lemma} in the rest of this section.

\qquad In order to show Lemma \ref{tournament_lemma}, we define an auxiliary process $(\mathtt{Rk}_{t})_{t\ge 0}$ on $\mathcal{T}_{x}$, $\mathtt{Rk}_{t}:\mathcal{V}\rightarrow \mathbb{Z}$ as follow: we first vertically set $\mathtt{Rk}_{0}(x)=\nbe_{t}(x)$ for all $t\ge 0$ and then horizontally extend to all vertices $\tilde{z}\in \mathcal{V}$ by 
\begin{equation}
	\mathtt{Rk}_{t}(\tilde{z}) = \mathtt{Rk}_{t}(x)+\int_{\vec{P}} \,d\tilde{X}_{t}  
\end{equation} 
for each $t\ge 0$. We are going to show that this new process satisfies Lemma \ref{tournament_lemma} and in fact equivalent to the original tournament expansion. We begin with the following commutative diagram which summarizes our construction.
\vspace{0.1cm}
\begin{equation}\label{diagram}
	\xymatrix{   
	 \text{ranking $\mathtt{Rk}_{t}$ on $\mathcal{T}_{x}$}\ar@{.>}[r] & \text{ranking $\mathtt{Rk}_{t+1}$ on $\mathcal{T}_{x}$}\\
            \text{edge configuration $d\tilde{X}_{t}$ on $\mathcal{T}_{x}$} \ar@{.>}[r]\ar[u] & \text{edge configuration $d\tilde{X}_{t+1}$ on $\mathcal{T}_{x}$}\ar[u]\\       
	\text{lift of 3-coloring $\tilde{X}_{t}$ on $\mathcal{T}_{x}$}\ar[u] \ar[r] & \text{lift of 3-coloring $\tilde{X}_{t+1}$ on $\mathcal{T}_{x}$}\ar[u]\\
  \text{3-coloring $X_{t}$ on $G$} \ar[r] \ar[u] & \text{3-coloring on $X_{t+1}$ on $G$} \ar[u]
	}
\end{equation}

In the above diagram, the first two horizontal maps from bottom are the CCA (or GHM) transition map on $G$ and $\mathcal{T}_{x}$, respectively, whereas the other two horizontal maps are to be described. In fact, the top horizontal map turns out to be the transition map for tournament processes. Since by definition $\mathtt{Rk}_{0}=\mathtt{rk}_{0}$, this will show the aforementioned equivalence. 

\qquad The induced dynamics $(d\tilde{X}_{0})_{t\ge 0}$ can be interpreted as a traditional comparison process called \textit{embedded particle system}, which consists of branching and annihilating particles on the edges of $\mathcal{T}_{x}$. This allows us to view propagation of excitation (or information flow) as moving and branching edge particles. Without loss of generality, we consider a CCA (or GHM) trajectory $(X_{t})_{t\ge 0}$ on a tree $T=(V,E)$ starting from an initial 3-coloring $X_{0}$. We view $\xi_{0}:=-dX_{0}$ as an initial \textit{edge configuration} where on each adjacent pair $(u,v)\in \overline{E}$ with $dX_{t}(u,v)=-1$, we place a single particle heading to $v$. We denote a particle on the edge $uv$ heading to $v$ by $u\rightarrow v$. The orientation of particles is so that 
\begin{equation}\label{particle_excitation_relation}
\text{$\exists$ a particle $v\rightarrow u$ at time $0$} \Longleftrightarrow   \text{$dX_{0}(v,u)=-1$} \Longleftrightarrow\text{$v$ excites $u$ at time $0$} 
\end{equation}

We let the edge particles evolve in discrete time in parallel, by applying the following rules successively for each transition $\xi_{t}\mapsto \xi_{t+1}$:

\begin{description}[noitemsep]
	\item{1.} (branching) If there is a particle $u\rightarrow v$ at time $t$, it branches into each incident edge $vw$ without a particle $v\leftarrow w$ at time $t$, and becomes a \textit{potential particle} $v\rightarrow w$ at time $t+1/3$.  
	\item{2.} (annihilation) If an edge has at least two potential particles at time $t+1/3$ with the opposite direction, then all potential particles on that edge get annihilated from the system at time $t+2/3$. 
	\item{3.} (coalescence) All remaining potential particles on each edge at time $t+2/3$, which have the same direction, coalesce into a single particle with the same direction at time $t+1$. Then $\xi_{t+1}$ counts the particles at time $t+1$. 
\end{description}

\begin{figure*}[h]
     \centering
          \includegraphics[width=0.95 \linewidth]{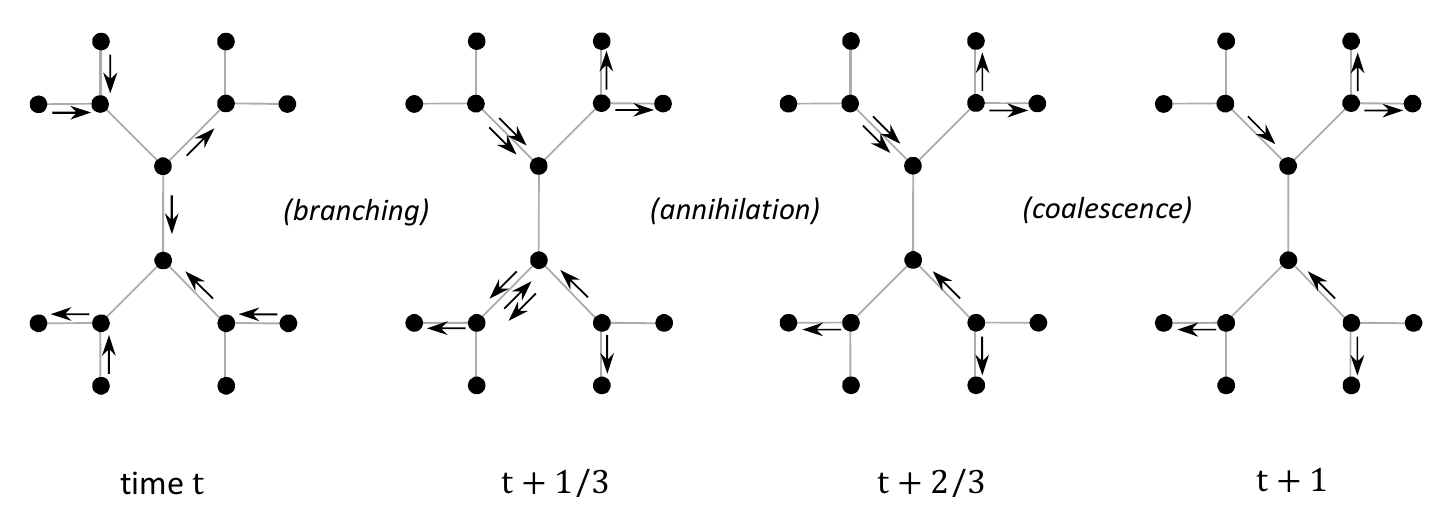}
          \caption{ Three-step edge particle evolution rule.
          }
          \label{aqps}
\end{figure*}
\vspace{-0.3cm}
We call this process $(\xi_{t})_{t\ge 0}$ the \textit{particle system expansion} of $(X_{t})_{t\ge 0}$. In the next proposition we show that this particle evolution rule is compatible with the actual CCA and GHM dynamics.  

\begin{customprop}{3.3}\label{particle_expansion}
	Let $T=(V,E)$ be a tree, $X_{0}:V\rightarrow \mathbb{Z}_{3}$ a 3-coloring, $(X_{t})_{t\ge 0}$ the CCA (or GHM) trajectory starting from $X_{0}$, and $(\xi_{t})_{t\ge 0}$ be its the particle system expansion. Then for each $t\ge 0$, we have the followings:
\begin{description}[noitemsep]
	\item{(i)} $dX_{t} = -\xi_{t}$.
	\item{(ii)} $(dX_{t+1}-dX_{t})(y,z)  = \indicator(\text{$z$ is excited at time $t$})- \indicator(\text{$y$ is excited at time $t$})$.	 
\end{description}
\end{customprop}

\begin{proof}
	First we show (i) implies  (ii). We show the following equivalent statement is implied from (i): 
\begin{equation}\label{flux_eq_base}
	\xi_{t+1}(y,z) - \xi_{t}(y,z) = \indicator(\text{$y$ is excited at time $t$})- \indicator(\text{$z$ is excited at time $t$}).
\end{equation}
Note that (i) says the relation (\ref{particle_excitation_relation}) continues to hold for all times $t\ge 0$. Hence $x$ is excited at time $t$ iff it has an incoming particle at that time. Suppose first that $\xi_{t}(y,z)=1$. By (i), $dX_{t}(y,z)=-1$, so $z$ is excited at time $t$. So it suffices to show that the first terms in (\ref{flux_eq_base}) agree. Indeed, since $\xi_{t}(y,z)=1$, all other potential incoming particles to $z$ at time $t$ get annihilated. Thus $\xi_{t+1}(y,z)=1$ iff $y$ has an incoming particle at time $t$ iff $y$ is excited at time $t$. Since both sides of (\ref{flux_eq_base}) are anti-symmetric, this deals with the opposite case $\xi_{t}(y,z)=-1$. Lastly, suppose $\xi_{t}(y,z)=0$. Then $\xi_{t+1}(y,z)=1$ if $y$ has an incoming particle at time $t$ but $z$ does not, $\xi_{t+1}(y,z)=-1$ if $z$ has an incoming particle at time $t$ but $y$ does not, and $\xi_{t+1}(y,z)=0$ if $y$ and $z$ are both excited or both not at time $t$. Hence (i) implies (ii).

\qquad It remains to prove (i). We prove the assertion for GHM dynamics, and a similar argument applies for CCA. The assertion holds for $t=0$ by definition. Suppose the assertion holds for $t\ge 1$. This means that at times $s\le t$, particles always point from color 1 to color 0. Now fix an edge $x_{1}x_{2}\in E$. First suppose $X_{t}(x_{1})=1$ and $X_{t}(x_{2})$=0. Then $X_{t+1}(x_{1})=2$ and $X_{t+1}(x_{2})=1$, so $dX_{t+1}(x_{1},x_{2})=0$. Hence we wish to show $\xi_{t+1}(x_{1},x_{2})=0$. In this case $dX_{t}(x_{1},x_{2})=-1$, and by the induction hypothesis, $\xi_{t}(x_{1},x_{2})=1$. Since particles are always from color 1 to 0 at times $\le t$, $x_{1}$ has no incoming particle at time $t$. Also, the particle $x_{1}\rightarrow x_{2}$ at time $t$ would annihilate all possible incoming particles to $x_{2}$ at time $t$. Thus $\xi_{t+1}(x_{1},x_{2})=0$ as desired. Second, suppose that both $X_{t}(x_{1})$ and $X_{t}(x_{2})$ are from $\{1,2\}$. Then $X_{t+1}(x_{i})=X_{t}(x_{i})+1 \mod 3$, so $X_{t+1}(x_{i})\ne 1$ for $i=1,2$, and $dX_{t+1}(x_{1},x_{2})=0$. Hence we also want to show $\xi_{t+1}(x_{1},x_{2})=0$. But this is clear since both $x_{i}$'s do not have incoming particles by the induction hypothesis. 

\qquad Lastly, suppose $X_{t}(x_{1})=X_{t}(x_{2})=0$. By the induction hypothesis, each $x_{i}$ has an incoming particle at time $t$ iff it has a neighbor of color 1 at time $t$. If neither  of them are excited at time $t$, then they both keep color 0 at time $t+1$, so $dX_{t+1}(x_{1},x_{2})=0$, which agrees with the particle dynamics since there are no incoming particles in either side. If both of them are excited, then $X_{t+1}(x_{1})=X_{t+1}(x_{2})=1$ so $dX_{t+1}(x_{1},x_{2})=0$. Also $\xi_{t+1}(x_{1},x_{2})=0$, as opposing particles come through both $x_{i}$'s and annihilate, so there is no remaining particles on the edge $x_{1}x_{2}$ at time $t+1$. Otherwise, by symmetry we may assume that $X_{t+1}(x_{1})=1$ and $X_{t+1}(x_{2})=0$, so $dX_{t+1}(x_{1},x_{2})=-1$. In terms of particles, $x_{1}$ has an incoming particle at time $t$ but $x_{2}$ does not, so there is a particle $x_{1}\rightarrow x_{2}$ at time $t+1$. This shows the assertion.    
\end{proof}

\qquad Next, we extend Proposition \ref{particle_expansion} (ii) to arbitrary cases. This gives the desired property of the tournament expansion  that a vertex increments its rank by 1 if and only if it gets excited in the original dynamics. 

\begin{customprop}{3.4}\label{excitation}
	Let $T=(V,E)$ be a tree and $(X_{t})_{t\ge 0}$, $(dX_{t})_{t\ge 0}$, and $(\mathtt{Rk}_{t})_{t\ge 0}$ be the four processes as before. Then for any $y,z\in V$ and $t\ge 0$, we have
\begin{equation}\label{flux}
	\int_{\vec{P}} \,dX_{t+1} - \int_{\vec{P}} \,dX_{t} = \indicator(\text{$z$ is excited at time $t$}) - \indicator(\text{$y$ is excited at time $t$})
\end{equation}	
where $\vec{P}$ is the unique directed path from $y$ to $z$ in $T$. Furthermore, for each $z\in V$ and $t\ge 0$, we have  
\begin{equation}
	\mathtt{Rk}_{t+1}(z)-\mathtt{Rk}_{t}(z) = \indicator(\text{$z$ is excited at time $t$}).
\end{equation}
\end{customprop}

\begin{proof} The second part of the assertion follows immediately from the first part, definition of ranking of the base point $x$, and the following relation coming from the definition: 
\begin{eqnarray}
	\mathtt{Rk}_{t+1}(z)-\mathtt{Rk}_{t}(z) = [\mathtt{Rk}_{t+1}(x)-\mathtt{Rk}_{t}(x)]+ \left[ \int_{\vec{P}} dX_{t+1} - \int_{\vec{P}} dX_{t} \right].
\end{eqnarray}
To show the first part, label the vertices of $\vec{P}$ by $y=x_{0},x_{1},\cdots,x_{k}=z$ such that $x_{i}x_{i+1}\in E$. Note that Proposition \ref{particle_expansion} (ii) gives 
\begin{eqnarray*}
	\sum_{i=0}^{k-1}(dX_{t+1}-dX_{t})(x_{i},x_{i+1}) &=& \sum_{i=0}^{k-1} \indicator(\text{$x_{i+1}$ is excited at time $t$})-\indicator(\text{$x_{i}$ is excited at time $t$})\\
	&=& \indicator(\text{$z$ is excited at time $t$})-\indicator(\text{$y$ is excited at time $t$})
\end{eqnarray*} 
which is equivalent to (\ref{flux}).  
\end{proof}

As a remark, this gives that contour integrals of $dX_{t}$ are time invariant. 

\textbf{Proof of Lemma \ref{tournament_lemma}.} The second assertion follows from the first since $\mathtt{rk}_{0}(x)=0$. To show the first part, note that by Proposition \ref{excitation} it suffices to show that 
\begin{equation}
	\mathtt{Rk}_{t}=\mathtt{rk}_{t} \quad \forall t\ge 0.
\end{equation}
It holds for $t=0$ by definition. Hence it suffices to show that each transition $\mathtt{Rk}_{t}\mapsto \mathtt{Rk}_{t+1}$ follows the transition map for tournament process. Indeed, by Proposition \ref{excitation} and construction, we have
\begin{eqnarray}
	\mathtt{Rk}_{t+1}(\tilde{z})-\mathtt{Rk}_{t}(\tilde{z}) 
	&=& \indicator(\text{$\tilde{z}$ is excited at time $t$}) \\
	&=&\indicator\left( \text{$\exists \tilde{y}\in N(\tilde{z})$ such that $d\tilde{X}_{t}(\tilde{y},\tilde{z})=-1$}\right)\\
	&=&\indicator\left( \text{$\exists y\in N(z)$ such that $\mathtt{Rk}_{t}(y)=\mathtt{Rk}_{t}(z)+1$} \right)
\end{eqnarray}
for any $\tilde{z}\in \mathcal{V}$ and $t\ge 0$. This shows the assertion. $\blacksquare$

\section{On finite graphs and the Erd\"os-R\'enyi random graph} 

We begin with a lemma that provides a lower bound for the activity in terms of contour integrals with respect to $dX_0$.
\begin{customlemma}{4.1}\label{lemma_irrotational}
	Let $G=(V,E)$ be arbitrary connected graph (not necessarily finite) with a cycle $\vec{C}$. Suppose $X_{0}:V\rightarrow \mathbb{Z}_{3}$ is an initial 3-coloring on $G$ such that the contour integral of $dX_{0}$ on $\vec{C}$ is nonzero. Then we have  
	\begin{equation}
		\alpha \ge \frac{1}{|V(\vec{C})|} \left| \oint_{\vec{C}} \,dX_{0} \right|. 
	\end{equation}     
\end{customlemma}

\begin{proof}
	 Without loss of generality, we may assume that $\oint_{\vec{C}}\,dX_{0}>0$. Fix $x\in V$, and let $(\mathtt{rk}_{t})_{t\ge 0}$ be the tournament expansion of $(X_{t})_{t\ge 0}$ at base point $x\in V$. By definition, Lemma \ref{tournament_lemma}, and Proposition \ref{tournamentdynamics}, we have 
	\begin{equation}
	\nbe_{t}(x) = \mathtt{rk}_{t}(x) = M_{t}(x).
	\end{equation}
	On the other hand, let $\vec{P}$ be any finite directed path in $G$ from $x$ to some vertex of $\vec{C}$, and let $\vec{W}_{n}$ be the walk starting from $x$, traversing $\vec{P}$, and then concatenating $\vec{C}$ $n$ times. Let $p=|E(\vec{P})|$ and $c=|V(\vec{C})|$. Then one has
	\begin{equation}
	M_{p+cn}(x)\ge \int_{\vec{W}_{n}}\,dX_{0} = \int_{\vec{P}}\,dX_{0}+n\int_{\vec{C}} \,dX_{0} 
	\end{equation}  
	so by taking $n\rightarrow \infty$, we obtain 
	\begin{equation}
	\alpha=\limsup_{t\rightarrow \infty} \frac{M_{t}(x)}{t} \ge \frac{1}{|V(\vec{C})|} \oint_{\vec{C}} dX_{0}
	\end{equation}
	as desired.
\end{proof}

\textbf{Proof of Theorem \ref{mainthm1}.} It suffices to show part (i) and (ii) of the assertion. First we show part (i). Suppose $X_{0}:V\rightarrow \mathbb{Z}_{3}$ is an initial 3-coloring such that $dX_{0}$ is irrotational. Fix an arbitrary node $x\in V$. Let $(\mathtt{rk}_{t})_{\ge 0}$ be the tournament expansion at $x$ of the CCA (or GHM) trajectory $(X_{t})_{t\ge 0}$. Let $D$ be the diameter of $G$. By Lemma \ref{key_lemma}, the assertion follows from 
\begin{equation}
	\lVert \mathtt{rk}_{0} \rVert_{\infty} \le D.
\end{equation}
Indeed, let $\vec{W}$ be any finite  walk in $G$ starting from $x$ to some vertex, say $z$. Let $\vec{P}$ be a shortest directed path from $x$ to $z$ in $G$. Then $|E(\vec{P})|\le D$, and by using the irrotationality of $dX_{0}$, we have 
\begin{equation}
	|\mathtt{rk}_{0}(\tilde{z})|=\left| \int_{\vec{W}}\,dX_{0} \right| = \left| \int_{\vec{P}}dX_{0}\right| \le |E(\vec{P})|\cdot \lVert dX_{0}\rVert_{\infty} \le D
\end{equation}
as desired. 

\qquad To show (ii), suppose that on some directed walk

 the contour integral of $dX_{0}$ does not vanish. We first note that $(X_{t})_{t\ge 0}$ converges to a periodic limit cycle since $G$ is finite and the dynamics is deterministic. Hence the limit of $E_{t}(x)/t$ as $t\rightarrow \infty$ exists. Recall that $\mathtt{rk}_{t}(x)=\nbe_{t}(x)=M_{t}(x)$ by definition, Proposition  \ref{tournamentdynamics}, Lemma \ref{tournament_lemma}. Hence by Lemma \ref{lemma_irrotational}, it suffices to show that 
\begin{equation}\label{M_t_sup}
	\lim_{t\rightarrow \infty} \frac{M_{t}(x)}{t} \le  \sup_{\vec{C}}\frac{1}{|V(\vec{C})|} \oint_{\vec{C}} dX_{0} .
\end{equation}
where the supremum runs over all directed cycles in $G$. 

\qquad  Choose $\vec{C}$ such that the right hand side of (\ref{M_t_sup}) attains its supremum. Let $\vec{W}_{t}$ be any walk of length $t$ in $G$ starting from $x$. Since $G$ is finite, $\vec{W}_{t}$ may have lots of self-intersections for $t\gg |V|$. Let $\vec{C}_{1},\cdots,\vec{C}_{n_{t}}$ be the sequence of cycles arising in $\vec{W}_{t}$, in chronological order as one traverses it. Note that the number of edges in $\vec{W}_{t}$ that are not used by $\vec{C}_{i}$'s is at most $|V|$. Hence by using the choice of $\vec{C}$, we get
\begin{eqnarray}
	\int_{\vec{W}_{t}}\,dX_{0} &\le&  |V|+\sum_{i=1}^{n_{t}} \left|  \int_{\vec{C}_{i}} \,dX_{0} \right|\\
	&\le & |V|+\sum_{i=1}^{n_{t}} \frac{|V(\vec{C}_{i})|}{|V(\vec{C})|} \left|  \int_{\vec{C}_{i}} \,dX_{0} \right| \\
	&\le & |V| + \frac{t}{|V(\vec{C})|} \int_{\vec{C}}\,dX_{0}.
\end{eqnarray}
Since $\vec{W}_{t}$ was arbitrary with length $t$, this yields 
\begin{equation}
	M_{t} \le |V|+\frac{t}{|V(\vec{C})|} \int_{\vec{C}}\,dX_{0}.
\end{equation}
Thus after dividing both sides by $t$ and letting $t\rightarrow \infty$, one has (\ref{M_t_sup}) as desired. $\blacksquare$   \\

We now proceed to the proofs of Theorems~\ref{mainthm2_cycle} and~\ref{corollary_ER}.  In both theorems, recall that the initial coloring is uniform. 

\textbf{Proof of Theorem \ref{mainthm2_cycle}.} Let $G=(V,E)$ be a connected graph, not necessarily finite, containing at least one cycle. It suffices to show the second assertion. Let $\{e_{1},\cdots,e_{k}\}\subset E$ be a matching in $G$ and let $C_{1},\cdots,C_{k}$ be cycles in $G$ such that $e_{i}\in E(C_{i})$ for $1\le i \le k$. Let $\vec{C}_{i}$ denote a directed cycle on $C_{i}$ with any of the two orientations given. According to the claim, we have 
\begin{equation}\label{mainthm_2_eq1}
	\mathbb{P}(\text{$X_{t}$ synchronizes weakly}) \le \mathbb{P}\left( \oint_{\vec{C}_{i}} \, dX_{0}=0 \quad \forall 1\le i \le k \right).
\end{equation}
Suppose all vertices but the $2k$ vertices used by the matching have been colored by $X_{0}$. Let $x_{i}$ and $y_{i}$ be the two endpoints of $e_{i}$ for each $1\le i\le k$. We are going to show that one can always assign colors on each pair $(x_{i},y_{i})$ with at least probability $2/9$ in such a way that $C_{i}$ is singular. This means that each $C_{i}$ can be non-singular with probability $\ge c$, and this occurs independently for each $i$. Thus the right hand side of (\ref{mainthm_2_eq1}) is at most $(7/9)^{k}$, as desired.

\qquad Fix $1\le i\le k$, and recall the definition of $dX_{0}$ in CCA and GHM cases, which was given at the beginning of Section 2. Let $x'_{i}$ and $y_{i}'$ be the neighbors of $x_{i}$ and $y_{i}$ on $C_{i}$, which are distinct iff $C_{i}$ has more than three vertices. Let $\vec{P}_{i}$ and $\vec{Q}_{i}$ be the two directed paths from $x_{i}'$ to $y_{i}'$ on $C$, where $\vec{P}_{i}$ is the one that contains $x_{i}$. Clearly $C$ is singular iff $\int_{\vec{P}_{i}}\,dX_{0}=\int_{\vec{Q}_{i}}\,dX_{0}$. First suppose CCA dynamics. By symmetry, we may assume that $(X_{0}(x_{i}'),X_{0}(y_{i}'))=(0,0)$ or $(0,1)$. Suppose $(X_{0}(x_{i}'),X_{0}(y_{i}'))=(0,0)$. If $\int_{\vec{Q}_{i}}\,dX_{0}=0$, then we set $(X_{0}(x_{i}),X_{0}(y_{i}))=(1,2)$ or $(2,1)$ so that $\int_{\vec{P}_{i}}\,dX_{0}=\pm 2$. Otherwise, we set $X_{0}(x_{i})=X_{0}(y_{i})$ so that $\int_{\vec{P}_{i}}\,dX_{0}=0$. Second, suppose $(X_{0}(x_{i}'),X_{0}(y_{i}'))=(0,1)$. In this case $\int_{\vec{Q}_{i}}\,dX_{0}\equiv 1\mod 3$, so it is enough to make $\int_{\vec{P}_{i}}\,dX_{0}$ to have the opposite sign. We can choose $(X_{0}(x_{i}),X_{0}(y_{i}))=(2,2)$ or $(0,2)$ to make the integral over $\vec{P}_{i}$ equals to $2$, or choose from $(0,0)$ or $(1,1)$ to make it $-1$. This takes care of the CCA case.

\begin{figure*}[h]
     \centering
          \includegraphics[width=0.7 \linewidth]{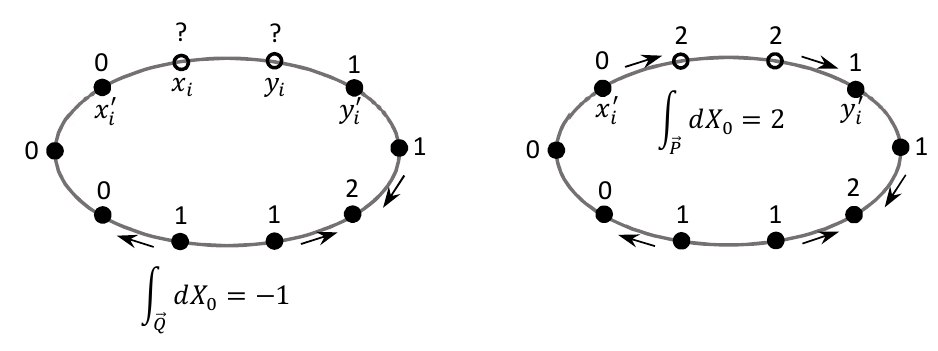}
          \caption{ Controlling singularity of a cycle using two adjacent vertices in case of CCA dynamics.
          }
          \label{aqps}
\end{figure*}

\qquad Now we suppose GHM dynamics. Since there is no symmetry in colors as in CCA, we need to consider all six cases of $X_{0}(x_{i}')\le X_{0}(y_{i}')$. We first suppose $\int_{\vec{Q}_{i}}\,dX_{0}\ne 0$. Note that $\int_{\vec{P}_{i}}\,dX_{0}= 0$ if $X_{0}(x_{i})=X_{0}(y_{i})=2$ regardless of the colors of their neighbors. Hence it suffices to find one additional choices for which $\int_{\vec{P}_{i}}=0$ in each of the six cases. We list them as tuples $(X_{0}(x_{i}'),X_{0}(x_{i}),X_{0}(y_{i}),X_{0}(y_{i}'))$; $(0,1,1,0)$, $(0,2,1,1)$, $(0,2,1,2)$, $(1,1,2,1)$, $(1,1,2,2)$, and $(2,1,1,2)$. It remains to check the case when $\int_{\vec{Q}_{i}}dX_{0}=0$. In this case, we find two tuples for which the integral over $\vec{P}$ is nonzero for each of the six cases. Namely, $(0,1,2,0)$, $(0,2,1,0)$, $(0,1,0,1)$, $(0,1,2,1)$, $(0,1,1,2)$, $(0,1,2,2)$, $(1,0,2,1)$, $(1,2,0,1)$, $(1,0,2,2)$, $(1,1,0,2)$, $(2,1,0,2)$, and $(2,0,1,2)$. This shows the assertion. $\blacksquare$

\qquad Now we turn our attention to Erd\"os-R\'enyi random graphs. By Theorem \ref{mainthm2_cycle} we know that the limiting behavior of $X_{t}$ is closely related to the existence of cycles and a matching that separates a large number of cycles. We begin preparing for the proof of Corollary $\ref{corollary_ER}$ with the following lemma.

\begin{customlemma}{4.1}\label{matching_ER}
	Let $G=G(m,\lambda/m)$ be the Erd\"os-R\'enyi random graph  for some $\lambda>0$. Let $M(G)$ be the size of the largest matching in $G$. Then for any $0<\delta<1$, there exists $m_{\delta}>0$ such that  
	\begin{equation}\label{matching_ER-eqn}
		\mathbb{P}\left( M(G) \le \frac{1}{3}(1-\delta)\mu m \right) \le e^{-\delta^{2}m/24}
	\end{equation}
for all $m> m_{\delta}$, where $\mu=1-e^{-\lambda/3}$.
\end{customlemma}

\begin{proof}
	Partition $V(G)=[m]$ into three sets $A,B$, and $C$ of equal size (up to rounding) and in particular $|B|=\lceil m/3 \rceil$. We describe an algorithm that finds a random matching in $G$ of some size $M_{m}$. We will then show that the random variable $M_{m}$ is large enough with sufficiently high probability to guarantee~\eqref{matching_ER-eqn}. 
	
	\qquad Label the vertices in $A$ as $x_{1},\cdots,x_{|A|}$. Initialize $A_{0}=A$, $B_{0}=B$, $C_{0}=C$, and $\mathcal{M}_{0}=\emptyset$. Define nested vertex sets $A_{i},B_{i}$ and $C_{i}$ and matching $\mathcal{M}_{i}$ recursively as follows:
	\begin{description}
		\item{(i)} If $x_{i}\in A_{i}$ it has a neighbor in $B_{i}$, say $y_{i}$, then set $\mathcal{M}_{i+1}=\mathcal{M}_{i}\cup \{x_{i}y_{i}\}$ and $A_{i+1}=A_{i}\setminus \{x_{1}\}$. Pick a vertex $z_{i}\in C_{i}$, and put $C_{i+1}=C_{i}\setminus \{z_{i}\}$ and $B_{i+1}=[B_{i}\setminus \{y_{i}\}]\cup \{z_{i}\}$.
		
		\item{(ii)} If $x_{i}\in A_{i}$ has no neighbor in $B_{i}$, then put $A_{i+1}=A_{i}\setminus \{x_{i}\}$, $B_{i+1}=B_{i}$, $C_{i+1}=C_{i}$, and $\mathcal{M}_{i+1}=\mathcal{M}_{i}$. 
	\end{description}	 
This process terminates after $|A|$ steps with matching $\mathcal{M}_{|A|}$. Let $M_{m}=|\mathcal{M}_{|A|}|$.

\qquad Now define random variables $X_{m,i}$, $1\le i \le |A|$, by 
\begin{equation}
	X_{m,i}=\indicator(\text{$x_{i}$ has a neighbor in $B_{i}$}).
\end{equation}
So $M_{m}=\sum_{i=1}^{|A|}X_{m,i}$. For each fixed $m\ge 1$, $X_{m,i}$ depends only on the size of $B_{i}$, which is constant and equals to $|B|$. Hence $X_{m,i}$'s are i.i.d. Bernoulli variables with $p_{m,i}=\mathbb{E}[X_{n,i}]=1-(1-\lambda/m)^{|B|}$. Since $|B|= \lceil m/3 \rceil$, we have  $p_{m,i}\rightarrow \mu:=1-e^{-\lambda/3}$ as $m\rightarrow \infty$. Now for a fixed $\delta>0$, there exists $m_{\delta}>0$ such that for all $m>m_{\delta}$, we have 
\begin{equation}
	\mathbb{P}\left[ \frac{M_{m}}{m} \le \frac{1}{3}(1-\delta)\mu \right] \le \mathbb{P}\left[ \frac{M_{m}}{|A|} \le (1-\delta/2)p_{m,i} \right] \le e^{-\delta^{2}|B|/8} \le e^{-\delta^{2} m/24},
\end{equation}
where the inequality in the middle is by Chernoff's bound for lower tail. This shows the assertion.
\end{proof}

\qquad The second ingredient to the proof of Theorem \ref{corollary_ER} is a large deviations estimate for the size of largest component in $G(n,\lambda/n)$ for $\lambda>1$, which we denote by $L(n,\lambda/n)$. It is well-known (see e.g. \cite{spencer1993nine}) that in this regime, $L(n,\lambda/n)/n$ converges in probability as $n$ tends to infinity, to the unique positive solution $\beta_{\lambda}$ to the equation  
\begin{equation}
\beta_{\lambda} = 1-e^{-\lambda \beta_{\lambda}}.
\end{equation}
We are interested in the large deviations estimate of the probability that the largest component in $G(n,\lambda/n)$ contains a fraction $x$ of vertices that is strictly less than the correct ratio $\beta_{\lambda}$. This is given by a reformulation of a more precise result by O'Connell (Lemma 3.2 in \cite{o1998some}). Namely, for each $0<x<\beta_{\lambda}$, we have 
\begin{equation}\label{ER_LDP}
\limsup_{n\rightarrow \infty} \frac{1}{n} \log \mathbb{P}[L(n,\lambda/n)<xn] = -A(x,\lambda)
\end{equation}
where
\begin{equation}
A(x,\lambda) = -x\log(1-e^{-\lambda x}) +x\log x + (1-x)\log (1-x) +\lambda x(1-x).
\end{equation}
We are now ready to give a proof of Theorem \ref{corollary_ER}.

\textbf{Proof of Theorem \ref{corollary_ER}.} We use standard facts about Erd\"os-R\'enyi random graph model, which may be found in many references including \cite{janson2011random}. If $p=o(1/n)$. Then we know that all components of $G$ are trees a.a.s., so (i) follows from Corollary \ref{corollary1}. To show (ii), suppose $p=\lambda/n$ for $\lambda\in (0,1)$. In this subcritical regime, every component in $G(n,p)$ is either a tree or contains a single cycle a.a.s. (see, e.g., Corollary 5.8 in Bollob\'as \cite{bollobas1998random}). Hence with Corollary \ref{corollary1}, we have  
\begin{eqnarray}
	&&\mathbb{P}(\text{$X_{t}$ synchronizes on each componenet of $G(n,p)$}) \\
	&& \qquad \qquad = \mathbb{P}(\text{$X_{t}$ synchronizes on all unicyclic components of $G(n,p)$}) + o(1). 
\end{eqnarray}

Now let $V_{ns}$ denote the total number of unicyclic components in $G(n,p)$ whose cycle is of length $s$ for $s\ge 3$. Let $\tau(x)=\sum_{s\ge 1} s^{s-1} x^{s}/s!$ be the exponential moment generating function of the sequence $\tau(k)=k^{k-1}$, $k\in \mathbb{N}$, which equals to the number of rooted trees on vertex set $[k]$. Pittel showed (see p.63 in \cite{pittel1988random}) that for each fixed $k\ge 3$, $V_{n3},\cdots,V_{nk}$ are asymtotically independent Poisson random variables with means $\nu_{3},\cdots,\nu_{k}$, where $\nu_{s}=\tau^{s}(\lambda e^{-\lambda})/2s$ for each $3\le s \le k$. Hence by Theorem \ref{mainthm1}, the assertion holds with 
\begin{equation}
	C(\lambda) = \prod_{k=3}^{\infty} \prod_{l=1}^{\infty} \frac{(\nu_{k})^{l}e^{-\nu_{k}}}{l!} \left[\mathbb{P}\left(\int_{C_{k}}\,dX_{0}=0\right)\right]^{l},
\end{equation}
where $C_{k}$ is a cycle of length $k$. It is not apparent however, from the formula above, that the constant $C(\lambda)$ is bounded away from both $0$ and $1$ for any $\lambda\in (0,1)$. To see this, we use Theorem 1 in Pittel \cite{pittel1988random}, which says 
\begin{equation}
	c(\lambda):=\lim_{n\rightarrow \infty}\mathbb{P}(\text{$G(n,p)$ has no cycle})= (1-\lambda)^{1/2}e^{\lambda/2+\lambda^{2}/4}\in (0,1).
\end{equation}
Thus by Theorem \ref{mainthm1} and \ref{mainthm2_cycle}, we have 
\begin{equation}
	0<c(\lambda) \le C(\lambda) \le 1-\frac{2}{9}(1-c(\lambda))<1.
\end{equation}

\begin{figure*}[h]
	\centering
	\vspace{-0.1cm}
	\includegraphics[width=0.45 \linewidth]{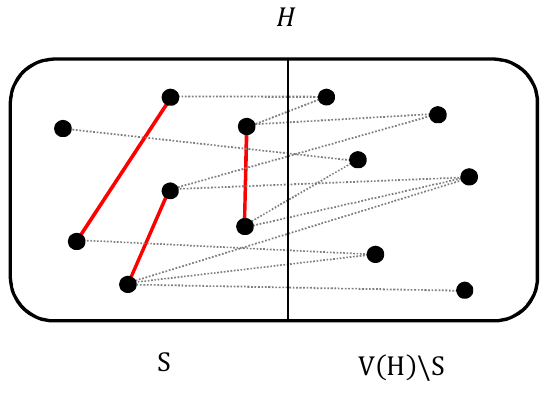}
	\vspace{-0.3cm}
	\caption{ Largest component $H\subset G$, a spanning tree $T$ of $H$ (dotted edges), a stable set $S\subset V(T)$, and a matching in $H[S]$ (red edges). 
	}
	\label{ER_pf_pic}
\end{figure*}

\qquad Next, suppose $p=\lambda/n$ for $\lambda>1$. We first introduce some notations. Let $H$ be the largest component of $G=G(n,\lambda/n)$. Choose a spanning tree $T$ of $H$, and let $S\subset V(H)$ such that no two vertices in $S$ have an edge in $T$ (i.e., a stable set in $T$). Since $T$ is bipartite, we may choose $S$ so that $|S|\ge \lceil |V(H)|/2 \rceil$. Finally, let $\mathtt{M}=M(H[S])$ be the size of largest matching in the induced subgraph $H[S]\subset H$ on the vertex set $S$. See Figure \ref{ER_pf_pic} for illustration.  

\qquad Define parameters $\lambda'=(\lambda+1)/2\in (1,\lambda)$, $\beta=\beta(\lambda)$ the unique solution of $\beta = 1-e^{-\lambda'\beta}$, and $c=c(\lambda):=\frac{\beta}{24}(1-e^{-(\lambda-1)\beta/24})$. Now define following events  
\begin{eqnarray*}
	A_{n} &=& \{ \text{$X_{t}$ synchronizes on the largest component of $G$}  \} \\
	B_{n} &=& \{ \text{The largest component $H\subseteq G$ has size $< \beta n/2$} \}\\
	C_{n} &=& B_{n}^{c} \cap \{\mathtt{M}< cn\}	
	\\
	D_{n} &=&  B_{n}^{c} \cap \{\mathtt{M}\ge cn\}	\cap \{\text{$X_{t}$ synchronizes on $H$}\}
\end{eqnarray*}
Then by partitioning, we may write 
\begin{equation}\label{partitioning}
	\mathbb{P}(A_{n})\le \mathbb{P}(B_{n})+\mathbb{P}(C_{n})+\mathbb{P}(D_{n}).
\end{equation}
  
\qquad Using (\ref{ER_LDP}), the first term in the right hand side of  (\ref{partitioning}) can be bounded by 
\begin{equation}
\mathbb{P}(B_{n}) \le \mathbb{P}[L(n,a/n)< \beta n/2] \le e^{-A(\beta/2,a)n}
\end{equation}
for large $n$. To bound the third term, observe that any matching $\{e_{1},\cdots,e_{k}\}$ in the induced subgraph $H[S]\subset H$ satisfies the hypothesis in Lemma \ref{matching_ER}. Indeed, for each $i$, let $C_{i}$ be the unique cycle which uses the edge $e_{i}$ and the unique path between its two endpoints in $T$ (see Figure \ref{ER_pf_pic}). Then for any $j\ne i$, $e_{j}\notin E(C_{i})$ since otherwise $e_{i}\in E(T)$, which contradicts the construction of $S$. Since on the event $D_{n}$ we have a matching in $G[S]$ of size at least $cn$, by Lemma \ref{corollary_ER} we have 
\begin{equation}
	\mathbb{P}(D_{n}) \le (7/9)^{cn}.
\end{equation}

\qquad Finally, let $a=\lambda'/n$ and $b=(\lambda-\lambda')/n$ and choose two independent random graphs $G(n,a)$ and $G(n,b)$. Their union has the law of $G(n,a+b-ab)$, which we may view as a subgraph of $G(n,\lambda/n)$ by a standard coupling. Denote by $G(S,b)$ the induced subgraph of $G(n,b)$ on the vertex set $S$. Note that 
\begin{equation}
G(\lceil xn/2 \rceil,b)\subseteq G(S,b) \subseteq G(S,\lambda/n) \quad \text{on the event $B_{n}^{c}$}
\end{equation}
after some possible relabeling of the vertices. Hence by Lemma \ref{matching_ER} with $\delta=1/2$ and $\mu=\lim_{n\rightarrow \infty} 1- e^{-\lceil \beta n/4 \rceil b/3} = 1-e^{-(\lambda-1)\beta /24}$, we have 
\begin{eqnarray}
	\mathbb{P}(C_{n}) &\le& \mathbb{P}\left[M( \text{$G(\lceil xn/2 \rceil,b)$ has no matching of size $\ge cn$} \right]\\
	& \le& e^{-\lceil \beta n/4 \rceil /96} \le e^{-\beta n/384}
\end{eqnarray} 
for all sufficiently large $n$. 

\qquad Combining the above estimates, we have $\mathbb{P}(A_{n})\le \exp[-Dn]$ for all sufficiently large $n$ where $D=D(\lambda)$ is given by 
\begin{equation}
	D(\lambda)=\frac{1}{384}\min\left\{ 192 A(\beta/2,(1+\lambda)/2),\,\beta,\,4\beta(1-e^{-(\lambda-1)\beta/24})\log(9/7)\right\}.
\end{equation}
$\blacksquare$

\section{On infinite trees}

\qquad In this section we prove Theorem \ref{mainthm2}. Throughout this section $\Gamma=(V,E)$ will denote an infinite but locally finite tree rooted at a vertex $0\in V$. None of our discussion depends on the choice of root. Let $X_{0}$ be a random 3-coloring on $\Gamma$ drawn from product measure $\mathbb{P}$ with marginal density $\mathbb{P}(X_{0}(\sigma)=i)=p_{i}$ for $i\in \mathbb{Z}_{3}$. Recall the definition of the $\Gamma$-indexed walk associated to the CCA or GHM dynamics $(X_{t})_{t\ge 0}$. Theorem \ref{mainthm2} gives quantitative estimates on the activity $\alpha$. The first step is to establish part (i) of the assertion, which states that 
\begin{equation}
	\alpha=v_{c} \qquad \text{a.s.}.
\end{equation}
According to Lemma \ref{key_lemma},  we have 
\begin{equation}\label{mainthm2_1_claim}
	\nbe_{t}(0) = \max_{|\sigma|\le t} S_{\sigma} \qquad \text{a.s. }
\end{equation}
for all $t\ge 0$. This implies 
\begin{equation}
	\alpha = \limsup_{t\rightarrow \infty}\frac{1}{t}\max_{|\sigma|\le t} S_{\sigma} \qquad \text{a.s.}
\end{equation}
Hence part (i) of Theorem \ref{mainthm2} follows from (\ref{mainthm2_1_claim}) and the following observation:

\begin{customprop}{5.1}\label{cloudspeed_nohistroy}
Almost surely, we have 
	\begin{equation}
		v_{c}=\limsup_{t\rightarrow \infty} \frac{1}{t} \max_{|\sigma|=t} S_{\sigma} = \limsup_{t\rightarrow \infty} \frac{1}{t} \max_{|\sigma|\le t} S_{\sigma}=\alpha.
	\end{equation}  
\end{customprop}

\begin{proof}
It suffices to show that 
\begin{equation}\label{v_c_pf}
	\limsup_{t\rightarrow\infty} \frac{1}{t}\max_{|\sigma|\le t} S_{\sigma}\le v_{c}.
\end{equation}
We may assume that $v_{c}< \infty$. Fix any $\epsilon>0$. Then there exists $T=T(\epsilon)>0$ such that 
\begin{equation}
	\max_{|\sigma|=t} S_{\sigma}\le (v_{c}+\epsilon)t \qquad \forall t>T.
\end{equation}
Let $M=\max_{|\sigma|\le T} S_{\sigma}$, and choose large $T'>0$ so that $T'>T$ and $M<(v_{c}+\epsilon)T'$. Then for all $t>T'$, we have 
\begin{eqnarray}
	\max_{|\sigma|\le t} S_{\sigma} &\le &\max_{|\sigma|\le T'} S_{\sigma} + \max_{T'< |\sigma|\le t} S_{\sigma}\\
	&\le & M+(v_{c}+\epsilon)t,
\end{eqnarray}     
which implies (\ref{v_c_pf}) as desired.  
\end{proof}

\qquad While Theorem \ref{mainthm2} (ii)-(v) gives an estimation on the cloud speed of the associated tree indexed walk, here we prove the assertion for general $\Gamma$-indexed walks with 1-correlated increments. Namely, let $Y_{\sigma}$ be a real-valued random variable attached to the edge $(\sigma^{-},\sigma)$, and assume throughout this section that they are identically distributed, centered, and that the moment generating function $\mathbb{E}[e^{tY_{\sigma}}]<\infty$ for some $t>0$. The $\Gamma$-index random walk $\{S_{\sigma}\}_{\sigma\in V}$ with increments $\{Y_{\sigma}\}_{\sigma\in V}$ is defined similarly:  
\begin{equation}
S_{\sigma} = \sum_{0\le \tau \le \sigma} Y_{\sigma} \qquad \forall \sigma\in V.
\end{equation}
Benjamini and Peres \cite{benjamini1994tree} obtained sharp upper and lower bounds on $v_{c}$ of a $\Gamma$-indexed walk for the case when $\Gamma$ with minimum degree $\ge 2$ and the increments are i.i.d. In the present work we generalize their result twofold: general infinite trees without degree constraint, and 1-correlated increments. Our arguments, which is largely based on Benjamini and Peres' original proof, easily carries over to any finite range correlation.

\qquad The analysis of cloud speed is essentially based on the large deviations principle for the random walk obtained by restricting the $\Gamma$-indexed walk on a single ray. Let $\gamma$ be an infinite ray in $\Gamma$ starting from $0$, and let $0=\sigma_{0},\sigma_{1},\sigma_{2},\cdots$ be the successive vertices on it. By the 1-correlation, $\{Y_{\sigma_{i}}\}_{i\ge 0}$ is a Markov chain, or more generally, a functional $g:\mathfrak{X}\rightarrow \mathbb{R}$ of some underlying Markov chain $\{\mathcal{X}_{i}\}_{i}$ on state space $\mathfrak{X}$. In our special case of CCA and GHM increments, we may take $\mathfrak{X}=(\mathbb{Z}_{3})^{2}$, $\mathcal{X}_{i}=(X_{0}(\sigma_{i}),X_{0}(\sigma_{i+1}))$, and $g$ determined by  
\begin{equation}
	Y_{\sigma_{i}}=dX_{0}(\sigma_{i},\sigma_{i+1}) = g(X_{0}(\sigma_{i}),X_{0}(\sigma_{i+1})).
\end{equation}
Furthermore, since $X_{0}$ is drawn from the uniform product measure $\mathbb{P}$ on $(\mathbb{Z}_{3})^{V}$, the product measure is the unique invariant measure for our chain and transition probabilities are given by 
\begin{equation}
	\pi[(a,b),(c,d)] = \begin{cases}
		1/3 & \text{if $b=c$}\\
		0 & \text{otherwise}.	
	\end{cases}
\end{equation}
Denote $\tilde{S}_{n}=\sum_{i=0}^{n-1}Y_{\sigma_{i}}$ with $\tilde{S}_{0}=0$. The large deviations principle for the measures associated to $\tilde{S}_{n}/n$'s deals with the probability that this sample average deviates from a typical behavior. When $Y_{\sigma}$'s are i.i.d., then this typical behavior is dictated by the strong law of large numbers, whereas ergodicity takes place in the Markovian context.  

\qquad To make a statement in the most convenient form for our purpose, let $\Lambda(t)=\log \mathbb{E}[e^{tY_{\sigma}}]$ be the logarithmic moment generating function of $Y_{\sigma}$ and $\Lambda^{*}$ be its \textit{Legendre transform} 
\begin{equation}
	\Lambda^{*}(u) = \sup_{t\in \mathbb{R}} [ut - \Lambda(t)].
\end{equation}
Note that this function may take $\infty$. When $Y_{\sigma}$'s are independent and satisfy the moment condition we gave at the beginning of this section, Cram\'er's theorem on large deviations  asserts that 
\begin{equation}\label{cramer}
	\lim_{n\rightarrow \infty} \frac{1}{n} \log \mathbb{P}[\tilde{S}_{n}\ge nu] = -\Lambda^{*}(u). 
\end{equation}
This formula is still valid when $\Lambda^{*}(u)$ takes the value $\infty$.

\qquad In order to extend this relation to Markovian setting, we may assume that the following (sufficient but not necessary) conditions are satisfied:
\begin{description}[noitemsep]
	\item{(a)} The state space $\mathfrak{X}$ is finite;
	\item{(b)} The Markov chain $\{\mathcal{X}_{i}\}_{i\ge 0}$ has a unique stationary distribution;
	\item{(c)} For each $t\in \mathbb{R}$, the principal eigenvalue $\lambda_{\pi}(tg)$ of a matrix $\pi_{tg}$ is positive, where $\pi_{tg}$ is the exponentially weighted transition matrix defined by 
\begin{equation}
	\pi_{tg}(x,y)=\pi(x,y)e^{tg(y)}.
\end{equation}
\end{description}
Then (\ref{cramer}) holds for the sums $\tilde{S}_{n}=\sum_{i=0}^{n-1}g(\mathcal{X}_{i})$ where $\Lambda$ is given by
\begin{eqnarray*}
	&& \Lambda(t):=\log \lambda_{\pi}(tg).
\end{eqnarray*}
We may refer to this relationship as Cram\'er's theorem for Markov chains. For references see~\cite{deuschel1989large, varadhan2008large}, and in particular Theorem 3.1.2 in~\cite{dembo2009large}.


\qquad We are now ready to state and prove a generalization of Proposition 4.1 in Benjamini and Peres \cite{benjamini1994tree}, from which the proof of Theorem \ref{mainthm2} will then follow.

\begin{customprop}{5.3}\label{cloudspeed_1correlated}
	Let $\Gamma=(V,E)$ be an infinite tree, and suppose the increments $\{Y_{\sigma}\}_{\sigma\in V}$ restricted on a ray starting from $0$ can be realized as a functional $g$ of some Markov chain $\{\mathcal{X}_{i}\}_{i\ge 0}$, which satisfies the conditions (a)-(c) above.  Then we have the followings:
	\begin{description}[noitemsep]
		\item{(i)} $v_{c}>0$ if and only if $\,\rentropy(\Gamma)>0$.
		\item{(ii)} $v_{c}$ satisfies the following upper bound 
	\begin{equation}\label{vc_upper_entropy}
		\Lambda^{*}(v_{c})\le \rentropy(\Gamma).
	\end{equation}
		\item{(iii)} If $\log \br(\Gamma)=\entropy(\Gamma)$, then either $v_{c}=B$ or equality in  (\ref{vc_upper_entropy}) is achieved, where $B:=\sup\{v\,:\, \Lambda^{*}(v)<\infty\}$. 
	\item{(iv)}  For each $r>1$ such that $v_{c}/(1-r^{-1})<B$, where $B$ is defined in (iii), we have  
	\begin{equation}\label{vc_lowerbd_rcut}
		\frac{\entropy_{r}(\Gamma)}{r-1}\le \Lambda^{*}\left( \frac{v_{c}}{1-r^{-1}}\right).
\end{equation}		
In particular, if $\Gamma$ has no leaves, then 
	\begin{equation}
		 \Lambda\left(\entropy(\Gamma)/v_{c}\right)\le \entropy(\Gamma),
	\end{equation}	
	which is sharp in the sense that there exists a tree-indexed walk on some tree $\Gamma^{-}$ on which its cloud speed $v_{c}^{-}$ satisfies $\Lambda(\entropy(\Gamma)/(v_{c}^{-}-\epsilon))>\entropy(\Gamma)$ for each $\epsilon>0$.
	\end{description} 
\end{customprop}

\begin{proof}
	\begin{description}
	\item{(i)} Follows from (ii) and (iv).
	\item{(ii)} By Cram\'er's theorem $v_{c}\le B$. We first show that $v_{c}<B$ implies $\Lambda^{*}(v_{c})\le \entropy(\Gamma)$, and improve the later result. Suppose $0\le v_{c}<B$. Let $\tilde{S}_{n}=\sum_{i=0}^{n-1}Y_{\sigma_{i}}$ be the one dimensional random walk as before. Observe that 
	\begin{eqnarray}
		\sum_{n\ge 1}\sum_{|\sigma|=n}\mathbb{P}(S_{\sigma}\ge nv) &=& 
		\sum_{n\ge 1} A_{n} \mathbb{P}(\tilde{S}_{n}\ge nv)  \\
		&=& \sum_{n\ge 1} \exp\left[\log A_{n} \mathbb{P}(\tilde{S}_{n}\ge nv) \right]   \\
		&=& \sum_{n\ge 1} \exp\left[n\left( \frac{1}{n}\log A_{n}+\frac{1}{n} \log \mathbb{P}(\tilde{S}_{n}\ge nv) \right) \right].
	\end{eqnarray}
Thus if $\entropy(\Gamma)<\Lambda^{*}(v)$, Cram\'er's theorem for Markov chains yields that the expression in the parenthesis above is negative and bounded away from 0 for large $n$. Thus the left hand side is summable, so by Borel-Cantelli lemma we get $v_{c}< v$. This shows $\Lambda^{*}(v_{c})\le \entropy(\Gamma)$. 

\qquad To improve the upper bound, suppose $v_{c}<B$. Fix $\epsilon>0$ and we construct a subtree $\Gamma_{\epsilon}$ of $\Gamma$ as follows. At each level $n=1,2,\cdots$ successively, delete all vertices $\sigma$ in $\Gamma$ if they have no descendant at level $\ge (1+\epsilon)n$. Call this operation $\epsilon$-pruning. Let $v_{c}(\epsilon)$ be the cloud speed of the $\Gamma$-indexed walk restricted on $\Gamma_{\epsilon}$. Since we can recover the original tree from $\Gamma_{\epsilon}$ by attaching deleted parts at levels $n=1,2,\cdots$ successively and since they all have depth at most $\epsilon n$, we get 
\begin{equation}
	v_{c}-\epsilon B \le v_{c}(\epsilon) \le v_{c}. 
\end{equation}
Moreover, since $\Lambda^{*}$ is non-decreasing and according to the first part, we know 
\begin{equation}
	\Lambda^{*}(v_{c}-\epsilon B)\le \entropy(\Gamma_{\epsilon}).
\end{equation}
Then by continuity of $\Lambda^{*}$ at $v=v_{c}$, it suffices to show that 
\begin{equation}\label{vc_upper_claim}
 \limsup_{\epsilon\rightarrow 0} \entropy(\Gamma_{\epsilon})\le \rentropy(\Gamma).
\end{equation}
To this end, recall that for any $n\le m$, $A_{n,m}$ denotes the number of vertices in $\Gamma$ at level $n$ which have descendants at level $m$. Let $A_{n}(\epsilon)$ denote the number of vertices in $\Gamma_{\epsilon}$ at level $n$. Note that $A_{n}(\epsilon)$ equals to the number of vertices in $\Gamma$ at level $n$ which survives the $\epsilon$-pruning upto level $n-1$ \textit{and} with descendants at level $\floor{(1+\epsilon)n}$. This yields 
\begin{equation}
	A_{n}(\epsilon) \le A_{n,\floor{(1+\epsilon)n}} 
\end{equation}  
for all $n\ge 1$, which easily implies $\entropy(\Gamma_{\epsilon})\le \entropy_{1+\epsilon}(\Gamma)$. This gives (\ref{vc_upper_claim}), as desired.

\item{(iii)} Here we  follow Lyons' argument in \cite{lyons1990random} with a minor modification. Suppose $\log\br(\Gamma)=\rentropy(\Gamma)=\entropy(\Gamma)$. Fix arbitrary $v>0$ such that $\Lambda^{*}(v)< \log \br(\Gamma)<\infty$. We will show that $v\le v_{c}$. Then the assertion follows since $\Lambda^{*}$ is strictly increasing on $[0,B]$. Choose $0<\epsilon<\entropy(\Gamma)-\Lambda^{*}(v)$. Then by Cram\'er's theorem for Markov chains, we may choose a  large $k\ge 1$ for which 
\begin{equation}\label{percolation_prob}
	\mathbb{P}\left[\tilde{S}_{k-1} \ge (k-1)v+B\right] > e^{-(\Lambda^{*}(v)+\epsilon)k} > e^{-\entropy(\Gamma)k} = (\text{br}\Gamma)^{-k}.
\end{equation}

For this choice of $k$, define a tree $\Gamma^{k}=(V_{k},E_{k})$ from $\Gamma$ by $V_{k}=\{\sigma\in V\,:\, k|\,|\sigma|\}$ and $(\sigma\rightarrow\sigma')\in E_{k}$ iff the unique path between $\sigma$ and $\sigma'$ in $\Gamma$, which we denote by $\gamma(\sigma,\sigma')$, has length $k$. In words, $\Gamma^{k}$ describes how vertices of $\Gamma$ at levels multiple of $k$ are interconnected. It is easy to check that $\text{br}\Gamma^{k}=(\text{br}\Gamma)^{k}$. We are going to define a bond percolation process on $\Gamma^{k}$, by deleting each edge $(\sigma\rightarrow \sigma')$ from $\Gamma^{k}$ if, roughly speaking, the partial sum on the corresponding path in $\Gamma$ grows slowly. More precisely, form a random subgraph $\Gamma^{k}(\omega)$ of $\Gamma^{k}$ by deleting each edge $(\sigma\rightarrow \sigma')$ unless 
\begin{equation}\label{edge_selection}
	\sum_{\sigma^{+} < \tau < \sigma' } Y_{\tau} \ge (k-1)v+B
\end{equation}
where $\sigma^{+}$ denotes the unique descendant of $\sigma$ on the path $\gamma(\sigma,\sigma')$. Note that the partial sum above is over the path $\gamma(\sigma^{+},\sigma')=\gamma(\sigma,\sigma')-\sigma$, and these paths are vertex disjoint in $\Gamma$. Hence by the 1-correlation, each edge in $\Gamma^{k}$ is selected independently. This defines a quasi-Bernoulli percolation process on $\Gamma^{k}$, which was introduced in \cite{lyons1989ising}. There it was shown that percolation occurs a.s. if $q_{k} \br(\Gamma^{k})>1$, where $q_{k}$ is the survival probability of each edge. Since this condition is satisfied by (\ref{percolation_prob}) and $\br(\Gamma^{k})=[\br(\Gamma)]^{k}$ in our case, $\Gamma^{k}(\omega)$ contains an infinite ray $\gamma$ a.s. This gives an infinite ray $\gamma_{0}$ in $\Gamma$ emanating from 0 such that all but finitely many vertices are in $\gamma$. Then $\limsup_{\sigma\in V(\gamma_{0})} S_{\sigma}/|\sigma|\ge v$, so this shows $v_{c}\ge v$ as desired.

	\item{(iv)} We follow the argument in Benjanini and Peres \cite{benjamini1994tree} with a minor modification. Suppose $\entropy_{r}(\Gamma)>0$ for some $r> 1$. For brevity we shall omit roundings in this proof. Fix $\epsilon>0$, and choose a strictly increasing $f:\mathbb{N}\rightarrow \mathbb{N}$ such that 
	\begin{equation}\label{condition_on_Ak}
		\lim_{n\rightarrow \infty} \frac{1}{f(n)} \log A_{f(n),rf(n)}>\entropy_{r}(\Gamma)-\epsilon. 
\end{equation}	 
For each $v>0$, define the following event 
	\begin{equation}
		\Omega_{n}(v) = \left\{\max_{|\sigma|=rf(n)} S_{\sigma} \le rf(n)v \right\}.
	\end{equation}		
For each vertex $\sigma$ at level $n$ which has a descendant, say $\sigma'$, at level $rf(n)$, let $\gamma(\sigma,\sigma')$ be the unique path between the two vertices. Note that there are at least $A_{f(n),rf(n)}$ of such paths, and they are mutually vertex disjoint. Denote by $M_{k}$ the random variable 
\begin{equation}
	M_{k}=\sum_{|\sigma|=k}\indicator(\text{$S_{\sigma}\ge -k\epsilon$   and $\sigma$ has a descendant at level $rk$}). 
\end{equation}
By conditioning on the values of $S_{\sigma}$ for $|\sigma|=f(n)$, we have 
\begin{equation}\label{path_splitting_at_level_k}
	\mathbb{P}\left[\Omega_{n}(v)\right] \le \mathbb{P}\left[M_{f(n)}\le A_{f(n),rf(n)}/2\right] + \mathbb{P}\left[\tilde{S}_{(r-1)f(n)} \le rf(n)v + f(n)\epsilon\right]^{ A_{f(n),rf(n)}/2}.
\end{equation}
The first term in the right hand side should tend to 0 as $k\rightarrow \infty$, since it is unlikely that for half of  $\sigma's$ at level $f(n)$ which has descendant at level $rf(n)$ we have $S_{\sigma}\le -f(n)\epsilon$, which is well below the correct mean $0$. Since
\begin{equation}
	\lim_{n\rightarrow \infty}\frac{\tilde{S}_{n}}{n} =0 \qquad \text{a.s., }
\end{equation}
we have 
\begin{equation}
	\frac{\mathbb{E}[M_{k}]}{A_{k,rk}} = \mathbb{P}(\tilde{S}_{k}\ge -k\epsilon) \rightarrow 1 \quad \text{as} \quad k\rightarrow \infty.
\end{equation}
But since $M_{k}\le A_{k,rk}$, this implies that 
\begin{equation}\label{odd_paths_estimate}
	\mathbb{P}(M_{k}\le A_{k,rk}/2) \rightarrow 0 \quad \text{as} \quad k\rightarrow \infty.
\end{equation}

\qquad Now we may assume $v_{c}<B(1-r^{-1})$ and pick $v\in (v_{c},B]$. The definition of cloud speed enforces that $\lim_{n\rightarrow \infty} \mathbb{P}(\Omega_{n}(v))=1$. Hence combining with (\ref{path_splitting_at_level_k}) and (\ref{odd_paths_estimate}), we get 
\begin{equation}
	\mathbb{P}\left[\tilde{S}_{(r-1)f(n)}\le f(n)(v+\epsilon)\right]^{A_{f(n),rf(n)}}\rightarrow 1 \qquad \text{as $n\rightarrow \infty$}
\end{equation}
Taking logarithm and and using the fact that $\log(1-x)\le -x$ for $|x|\ll 0$, this yields  
\begin{equation}\label{mainthm2_eq1}
	A_{f(n),rf(n)}\mathbb{P}\left[ \tilde{S}_{(r-1)f(n)}\ge f(n)(v+\epsilon)\right]\rightarrow 0 .
\end{equation}
By Cram\'er's theorem for Markov chains, if $(r-1)f(n)$ is sufficiently large, then 
\begin{equation}\label{cramer_ineq1}
	\mathbb{P}\left[ \tilde{S}_{(r-1)f(n)}\ge f(n)(v+\epsilon) \right] \ge \exp[-(r-1)f(n)(\Lambda^{*}(u)+\epsilon)], 
\end{equation}
where 
\begin{equation}\label{def_u}
	u = \frac{v+\epsilon}{1-r^{-1}}.	
\end{equation}
For large enough $n$, the terms in (\ref{mainthm2_eq1}) is less than 1. At the same time, by (\ref{condition_on_Ak}), we have $\log A_{f(n),rf(n)}\ge f(n)(\entropy(\Gamma)-\epsilon)$ for $n$ large. Thus (\ref{mainthm2_eq1}) and (\ref{cramer_ineq1}) gives 
\begin{equation}\label{mainthm2_ineq2}
	(\entropy_{r}(\Gamma)-\epsilon)\le (r-1)(\Lambda^{*}(u)+\epsilon).
\end{equation}
Then by the continuity of $\Lambda^{*}$ at $v_{c}<B$, letting $v\searrow v_{c}$ and $\epsilon\searrow 0$ establishes the first assertion. 

\qquad To show the second part, we may further assume that $\Gamma$ has no leaves. Then we have $\entropy_{r}(\Gamma)=\entropy(\Gamma)$ for all $r> 1$. We may write (\ref{mainthm2_ineq2}) as 
\begin{equation}\label{mainthm2_legendre}
	\frac{\entropy(\Gamma)-\epsilon}{v+\epsilon}u -\Lambda^{*}(u)\le \entropy(\Gamma)
\end{equation}
where $u$ is given by (\ref{def_u}). Since this is valued for all $r> 1$, the above inequality is valid for all $u>v+\epsilon$. Since it is also valid for $0<u\le v+\epsilon$ trivially, we may take supremum of left hand side of (\ref{mainthm2_legendre}) over all $u>0$, which makes it the Legendre transform of $\Lambda^*$. But since $\Lambda$ is convex, we have $\Lambda^{**}=\Lambda$ by the Fenchel-Moreu theorem \cite{arnol2013mathematical}. This gives 
\begin{equation}
	\Lambda\left( \frac{\entropy(\Gamma)-\epsilon}{v+\epsilon} \right) \le \entropy(\Gamma),
\end{equation}
and by left continuity of $\Lambda$, letting $v\searrow v_{c}$ and $\epsilon\searrow 0$ establishes second assertion.

\qquad Lastly, to show the lower bound is sharp, consider a tree $\Gamma^{-}=(V,E)$ where every vertex has only one descendant except one special vertex $\sigma_{j}$ at every level $j!$ for each $j\ge 1$, which has $e^{dj!}$ descendants. This is called the ``exploding tree'' in \cite{benjamini1994tree}, and in the reference it is shown that the cloud speed on this tree satisfies the assertion. This completes the proof. 
\end{description}		
\end{proof}

{\bf Proof of Theorem~\ref{mainthm2}.} Part (i) follows from Proposition~\ref{cloudspeed_nohistroy}.  By Proposition~\ref{cloudspeed_1correlated}, it remains to verify that $\Lambda$ in the cases of CCA and GHM agree with $\Lambda_{\text{CCA}}$ and $\Lambda_{\text{GHM}}$, as stated in Section~\ref{section:results}.  The corresponding exponentially weighted transition matrices $\pi_{tg}^{\text{CCA}}$ and $\pi_{tg}^{\text{CCA}}$ are 
\begin{equation}
\pi_{tg}^{\text{CCA}} = \begin{bmatrix}
p_{0}B_{1} & p_{1}e^{-t}B_{1} & p_{2}e^{t}B_{1} \\
p_{0}e^{t}B_{2} & p_{1}B_{2} & p_{2}e^{-t}B_{2} \\
p_{0}e^{-t}B_{3} & p_{1}e^{t}B_{3} & p_{2}B_{3} 
\end{bmatrix},
\quad 
\pi_{tg}^{\text{GHM}} = \begin{bmatrix}
p_{0}B_{1} & p_{1}e^{-t}B_{1} & p_{2}B_{1} \\
p_{0}e^{t}B_{2} & p_{1}B_{2} & p_{2}B_{2} \\
p_{0}B_{3} & p_{1}B_{3} & p_{2}B_{3} 
\end{bmatrix},
\end{equation}
where $B_{j}$ is the $3\times 3$ matrix with all zeros but ones in the $j^{\text{th}}$ column, for $1\le j \le 3$. Hence $\Lambda(t)$ is the logarithm of the largest positive root $x=x(t)$ in the corresponding characteristic polynomial. Elementary computations show that the resulting equations are as follows:
\begin{equation}
(\text{CCA})\qquad	x^{3}-x^{2}=p_{0}p_{1}p_{2}(e^{3t}+e^{-3t}-2)
\end{equation} 
\begin{equation}
(\text{GHM})\qquad x^{3}-x^{2}=p_{0}p_{1}p_{2}(e^{t}+e^{-t}-2)
\end{equation}
This shows that $\Lambda$ for CCA and GHM agree with $\Lambda_{\text{CCA}}$  and $\Lambda_{\text{GHM}}$, and that  $\Lambda_{\text{GHM}}(t)=\Lambda_{\text{CCA}}(t/3)$. $\blacksquare$\\

\begin{customexample}{5.5}\label{$d$-ary tree}
	
	($d$-ary tree) Fix an integer $d\ge 2$, and let $\Gamma=(V,E)$ be the $d$-ary tree. Suppose uniform color density for $X_{0}$. By part (iv) of Theorem \ref{mainthm2} and using the closed form expression for $\Lambda_{\text{CCA}}^{*}$, we compute the activity $\alpha$ of CCA dynamics on $\Gamma$ as 
	\begin{equation}
	\alpha = 0.86824163\cdots 
	\end{equation}
	for $d=2$ and $\alpha\equiv 1$ for $d\ge 3$. In case of GHM dynamics, part (iv) of Theorem \ref{mainthm2} gives 
	\begin{equation}
	\alpha = 0.28941386\cdots
	\end{equation}
	for $d=2$ and $\alpha\equiv 1/3$ for $d\ge 3$.   $\blacktriangle$\\
\end{customexample}

\begin{customexample}{5.6}\label{geometric tree}

(Geometric tree) Fix an integer $d\ge 2$, and let $\Gamma^{d}=(V,E)$ be obtained by a single infinite ray $\gamma$ emanating from root $0$, by attaching a $d$-ary tree of depth $ 4^{k}$ to each vertex $\sigma\in V(\gamma)$ at level $|\sigma|=4^{k}$. Consider CCA dynamics on $\Gamma$, so that $v_{c}\le 1$. Since $\Gamma^{d}$ has a single infinite branch, one has $\br(\Gamma^{d})=1$. It is also easy to check that $\entropy(\Gamma)=\frac{\log d}{2}$. Moreover, for any $r\in [1,2]$, $\Gamma^{d}$ has a $r$-cut $f$ defined by $f(k)=(2/r)4^{k}$.  This gives $\entropy_{r}(\Gamma)\ge (1-r/2)\log d$, so Proposition \ref{cloudspeed_1correlated} (iii) gives 
\begin{equation}
	\frac{1-r/2}{r-1}\log d \le \Lambda^{*}\left( \frac{v_{c}}{1-r^{-1}} \right). 
\end{equation}
Thus $v_{c}$ attains the maximum value 1 for large $d$, despite $\br(\Gamma)=1$. Also note that the optimal value of $r$ for the above lower bound is strictly between $1$ and $2$. This makes sense because if $r=1$, we have the largest $r$-entropy which matches the volume entropy, but there is not enough room for the maximum to grow. On the other hand, $r=2$ gives the longest time for the maximum to grow, but the corresponding $r$-entropy is the lowest.  In fact, in this example it is possible to calculate $v_{c}$ explicitly. By the second part of Proposition \ref{cloudspeed_1correlated} (ii), we know that the cloud speed $v_{c}'$ on a $d$-ary tree is given by $\Lambda^{*}(v_{c}')=\log d$. Observe that on $\Gamma^{d}$, for each interval $[4^{k},4^{k+1}]$ of levels, the maximum grows linearly with speed $v_{c}'$ only for the first third and stays the same for the later two thirds. Thus $v_{c}=v_{c}'/3$.  $\blacktriangle$
\end{customexample}

\section{Concluding remarks and open problems}

\subsection{CCA With more than three colors} 

\qquad Our key tool in characterizing the limiting behavior of 3-color CCA  is the tournament expansion, which enables us to unfold the dynamics of cyclically arranged colors into a monotone dynamics of linearly ordered ranks. An interesting and natural question is whether a similar technique might work
for $\kappa\ge 4$. In  this case, the particle representation is more complex \cite{fisch1990one}, and existence 
of a suitable tournament expansion remains unclear. Proving (or disproving) the following conjecture 
would be a key step in understanding these dynamics. 

\begin{customconjecture}{6.1}  
Fix any $\kappa\ge 4$. 
Let $\Gamma=(V,E)$ be a tree with  $\br(\Gamma)>1$, and let $X_{0}$ be a random $\kappa$-coloring of 
$V$ drawn from the uniform product measure. Then, almost surely, every vertex oscillates. 
\end{customconjecture}

Methods from the present paper yield some partial results. For example, 
it is not hard to prove the above conjecture for $\kappa=4$ and regular binary tree $\Gamma$. 
For very large $\kappa$, we expect that the activity is caused by large connected sets in which all edges have color differences $0$ or $\pm 1$. Once large enough, these regular droplets grow without bounds with high 
probability. This nucleation scenario seems to hold even for some trees with
branching number $1$ (see Figure 5), but, due to correlations in the
growth environment, techniques for proving its validity remain elusive.

\begin{figure*}[h]
	\centering
	\includegraphics[width=1 \linewidth]{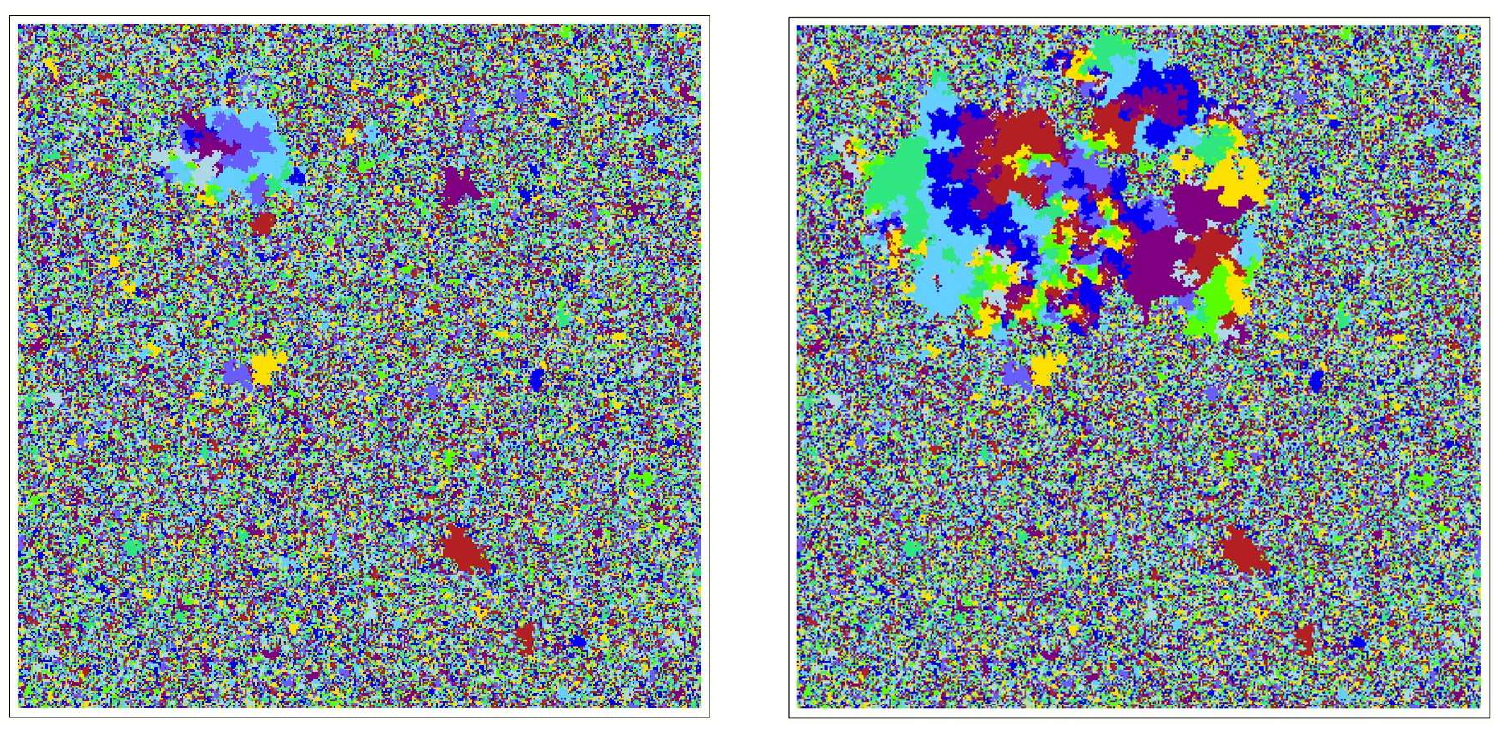}
	\caption{ Two snapshots of 9-color CCA dynamics on a uniform spanning tree of a $400\times 400$ torus, at times about 3,000 and 17,000. By the former time, almost the entire square fixates, except for a single droplet that continues to grow slowly until it takes over the available space.
	}
	\label{8CCA_snapshots}
\end{figure*}

\vskip0.3cm

\subsection{Activity for CCA and GHM on arbitary trees}

\qquad In Theorem \ref{mainthm2} (iv), we have seen that the activity for GHM is exactly a third of that for CCA for a given random 3-coloring $X_{0}$ with arbitrary color density, provided 
the tree is regular enough that $\log \br(\Gamma)=\entropy(\Gamma)$. In fact, it is not hard to see that an inequality holds for an arbitrary infinite tree $\Gamma$ in the special case of uniform density. Namely, let $\alpha_{\text{CCA}}$ and $\alpha_{\text{GHM}}$ be the activities of CCA and GHM dynamics on $\Gamma$, starting from a random 3-coloring $X_{0}$ with uniform density. Then we claim that 
\begin{equation}\label{claim_activity_ineq}
\alpha_{\text{CCA}}\le  3\alpha_{\text{GHM}}.
\end{equation}
Indeed, let $X_{0}^{(i)}$ for $i\in \mathbb{Z}_{3}$ be the random 3-coloring on $\Gamma$ obtained by color shifts from $X_{0}$, i.e., $X_{0}^{(i)}(x)=X_{0}(x)+i\mod 3$. Since $d_{\text{CCA}}X_{0} = \sum_{i=0}^{2} d_{\text{GHM}}X_{0}^{(i)}$, taking path integrals and maxima on both sides gives
\begin{equation}\label{claim_activity_pf}
\alpha_{\text{CCA}} \le \alpha_{\text{GHM}}^{(0)}+ \alpha_{\text{GHM}}^{(1)} + \alpha_{\text{GHM}}^{(2)},
\end{equation}
where $\alpha_{\text{GHM}}^{(i)}$ is the activity of GHM
with initial coloring $X_{0}^{(i)}$ for $i\in \mathbb{Z}_{3}$. As the initial color densities are uniform, $X_{0}^{(i)}$'s are identically distributed so that $\alpha_{\text{GHM}}^{(i)}$'s coincide,
and therefore (\ref{claim_activity_pf}) implies (\ref{claim_activity_ineq}).  

\qquad The message of Theorem \ref{mainthm2} (iv) is that, provided $\log \br(\Gamma)=\entropy(\Gamma)$, the correlation between the three GHMs in the above 
paragraph does not increase the cloud speed for the CCA, resulting in 
equality in (\ref{claim_activity_ineq}). To make sense of this result, observe from the percolation arguments in Section 5 that, on such regular trees, the cloud speed is essentially attained on a single ray. By contrast, on irregular trees the cloud speed might be determined 
by disconnected chunks. We thus pose by the following question.

\begin{customquestion}{6.2}
Let $\Gamma=(V,E)$ be an arbitrary infinite tree and let $X_{0}$ be a random $3$-coloring of $V$ drawn from the uniform product measure. Is it true that 
\begin{equation}
\alpha_{\text{CCA}} = 3\alpha_{\text{GHM}}?
\end{equation}
\end{customquestion}

One may attempt to resolve this issue by generalizing the characterization of the cloud speed $v_c$ from Theorem \ref{mainthm2} (iv) to arbitrary trees.  However, the cloud speed on irregular trees may not be characterized solely by a dimensionality of the underlying tree and step size distribution: there are two trees $\Gamma_{1}$, $\Gamma_{2}$ and two step size distributions $F_{1}$, $F_{2}$ such that $v_{c}(\Gamma_{1},F_{1})>v_{c}(\Gamma_{2},F_{2})$ but $v_{c}(\Gamma_{1},F_{2})<v_{c}(\Gamma_{2},F_{1})$ (see Remark in \cite{benjamini1994tree} following Proposition 4.3). Hence characterizing the cloud speed on irregular trees, which may be of an independent interest, seems to require novel techniques.

\section*{Acknowledgements} The authors are grateful to Russ Lyons for pointing out an error in an earlier draft. Janko Gravner was partially supported by the NSF grant DMS-1513340, 
Simons Foundation Award \#281309,
and the Republic of Slovenia's Ministry of Science
program P1-285. Hanbaek Lyu was partially supported by a departmental research fellowship. David Sivakoff was partially supported by NSF CDS\&E-MSS Award \#1418265.


\bibliographystyle{elsarticle-harv}   
\bibliography{mybib}  




\end{document}